\documentclass[11pt]{amsart}
\usepackage{amsmath,amssymb,amsthm}
\usepackage[margin=1.1in]{geometry}
\usepackage{xcolor}
\definecolor{linkblue}{RGB}{0,51,153}
\usepackage[colorlinks=true,linkcolor=linkblue,citecolor=linkblue,urlcolor=linkblue]{hyperref}
\newtheorem{theorem}{Theorem}[section]
\newtheorem{lemma}[theorem]{Lemma}
\newtheorem{proposition}[theorem]{Proposition}
\newtheorem{corollary}[theorem]{Corollary}
\theoremstyle{definition}

\newtheorem{example}[theorem]{Example}
\newtheorem{remark}[theorem]{Remark}
\newtheorem{problem}[theorem]{Problem}
\newcommand{\C}{\mathbb C}
\newcommand{\A}{\mathcal A}
\newcommand{\K}{\mathcal K}
\newcommand{\R}{\mathcal R}
\newcommand{\Par}{\mathcal P}
\newcommand{\gd}{\mathrm{gd}}
\newcommand{\Dz}{\mathrm D}
\newcommand{\ind}{\operatorname{ind}}
\newcommand{\rk}{\operatorname{rank}}
\newcommand{\tr}{\operatorname{tr}}
\newcommand{\im}{\operatorname{im}}
\newcommand{\adj}{\operatorname{adj}}
\newcommand{\gap}{\operatorname{gap}}
\title[Parallel sums of nilpotent elements]{Parallel sums of nilpotent elements:\\ exact realization and structured eigenvalue splitting}
\author{Yaashaa Golovanov}
\address{AMMOC, An International Math Circle in Honor of Vladimir Igorevich Arnold and Jerrold Eldon Marsden}
\email{golovanov@ammoc.org}
\date{November 25, 2026}
\begin{document}
\raggedbottom
\begin{abstract}
We prove that, over every infinite field and for $m\ge2$, the parallel sums
$N(N+C)^{-1}C$ of nilpotent $m\times m$ inputs with invertible sum
are exactly the matrices of nullity at least two. Every positive split
of the target nullity between the input kernels is realizable.
Over $\C$, we determine the paired derivative image and characterize
universal full matrix cancellation at nonzero residues.
For analytic input paths, nilpotence is required only at the base.
An adjugate identity forces the last two characteristic coefficients
to be $O(\varepsilon^2)$, yielding sharp bounds
$O(|\varepsilon|^{2/3})$ at nonzero three-dimensional nilpotent
residues and $O(|\varepsilon|^{1/(m-2)})$ at nilpotent residues
for $m\ge4$. At every base with residue of type $[m-1,1]$, $m\ge4$, we prove
a generic splitting law, also with fixed input sum. For $m\ge5$,
it has $m-2$ roots on the scale $|\varepsilon|^{1/(m-2)}$ and
two on the scale $|\varepsilon|^{1/2}$, separating the largest-root
and spectral-gap scales. A separate sharpness family keeps both
inputs nilpotent and two roots zero. At the nonzero nilpotent
three-dimensional boundary, a finite coefficient criterion computes the
splitting orders for each exceptional affine fixed-sum direction.
For each such path, the punctured spectral gap is at least
$c|\varepsilon|^3$ near zero for some $c>0$, and the exponent is
optimal over this class. Over a nonzero unital ring
with $2$ invertible, we construct two $3\times3$ matrix inputs,
each of nilpotency index exactly three, with invertible sum.
In the complex unital Banach case, the residue has spectrum
$\{0\}\cup\tfrac12\sigma(Q)$ for every coefficient $Q$.
With the common core allowed to vary,
a complex matrix is a finite-Drazin Hua output if and only if
its nullity is not one.
\end{abstract}
\maketitle

\section{Introduction}\label{sec:introduction}
Let $N$ and $C$ be nilpotent elements whose sum is invertible. The map
\begin{equation}\label{eq:parallel-map}
 \Par(N,C)=N(N+C)^{-1}C
\end{equation}
can have a much broader spectral range than either input. We call
$\nu=\Par(N,C)$ the \emph{residue}, or \emph{residual parallel sum},
of the pair. We ask which residues occur and how their eigenvalues
split when the inputs are perturbed.
The same parallel-sum expression occurs in the common-core reduction
of Hua inversion, where the complementary residue controls success
of the prescribed core inverse; Section~\ref{sec:corners} specifies
the finite-Drazin and Banach criteria.

Theorem~\ref{thm:complete-realization} answers the realization question
over every infinite field: the targets in $M_m$ are exactly the
matrices of rank at most $m-2$. It realizes every positive split of the
target nullity between the input kernels. In every realizing pair,
nullity two is equivalent to both inputs being regular nilpotent.
The proof works in characteristic two; Remark~\ref{rem:realization-F2}
gives the precise obstruction over $\mathbb F_2$ in dimension three.
Proposition~\ref{prop:index} gives the sharp nilpotent-output index
bound $m-1$ and explicit weighted-cycle witnesses.

For perturbations we work over $\C$. A paired path varies
$(a,c)$ with invertible sum, while its equivalent raw coordinates
are $(a,b)$ with $b=(a+c)^{-1}$ and output $F(a,b)=a-aba$.
Only the base inputs must be nilpotent. Fixed coupling means that
$a+c$ is constant; affine fixed coupling further specifies
$a=N+\varepsilon X$, $c=C-\varepsilon X$.
Section~\ref{sec:mechanism} makes these classes precise, including
their distinction from paths with selected inverse coordinates.
Proposition~\ref{prop:paired-derivative} determines the entire paired
derivative image, of codimension $2\dim\ker N\dim\ker C$.
The adjugate identity in Proposition~\ref{prop:adjugate} additionally
makes the last two characteristic coefficients $O(\varepsilon^2)$.
The matrix coefficient in the resulting first-order trace identity has
rank one at zero-primary type
$[h-1,1]$ with $h\ge3$, but rank two when $h=2$.
The ordered annihilation alternative is asserted only in the former
case (Proposition~\ref{prop:zero-primary}).

Section~\ref{sec:rates} develops the resulting splitting laws.
At a nonzero nilpotent residue in dimension three, every analytic
raw path has all roots $O(|\varepsilon|^{2/3})$; a first-derivative
coefficient distinguishes exact $2/3$ splitting from an
$O(|\varepsilon|)$ bound. The simultaneous normal form identifies
the two exceptional fixed-coupling hyperplanes. For $m\ge4$,
the universal bound at a nilpotent residue is
$O(|\varepsilon|^{1/(m-2)})$. Theorem~\ref{thm:extremal-splitting}
proves its generic attainment at every base of type $[m-1,1]$,
also within fixed coupling. For $m\ge5$, $m-2$ roots have that
scale and the remaining two have square-root scale. Formula
\eqref{eq:nilpotent-preserving-sharpness} separately proves sharpness
with both inputs remaining nilpotent and two roots remaining zero.
The section also treats the exact canonical cubics, the
two-zero-primary refinement and the rank-two square-zero boundary
under their respective hypotheses.
Proposition~\ref{prop:affine-coefficient-criterion} computes the
exceptional affine fixed-coupling rates from finitely many coefficients.
In this exceptional affine class, Corollary~\ref{cor:affine-sharp-gap}
gives the sharp cubic lower bound for the punctured gap and the sharp
order bound $5/2$ when no branch is identically zero.
Problem~\ref{prob:rate} retains the geometry of the resulting loci.

The two-thirds exponent itself is classical: Moro, Burke and
Overton~\cite[pp.~805 and~811]{MBO97} discuss the type-$[2,1]$ case.
Here the restrictions come from the nilpotent input map. Even affine
input paths produce nonlinear output paths, whose later terms matter
on the exceptional locus. The universal full matrix cancellation
criterion for unrestricted paired or raw directions in
Proposition~\ref{prop:full-cancellation} is also stronger than the
scalar trace cancellation: it is equivalent to one of two ordered
annihilation patterns and forces a square-zero output. Its necessary
dimension bound $m\ge3\rk\nu$ is attainable, while dimension four
with nilpotent base inputs and square-zero rank-two output instead has
generic square-root splitting.

The realization problem differs from the classical questions about
sums and products of nilpotent matrices. Over a field, products of
two nilpotent square matrices are precisely the singular matrices,
apart from the nonzero nilpotent matrices of order two; see
\cite[Theorem~4]{Hattingh16} for a corrected proof. Sums of two
nilpotent matrices are precisely the trace-zero matrices that are
zero or nonscalar; see \cite[Proposition~3]{BreazCalugareanu17}.
We instead prescribe the parallel sum and require the input sum
to be invertible. The cyclic characteristic-polynomial prescription
used in our proof is credited at its point of use.

Parallel addition originates in the network matrix theory of
Anderson and Duffin~\cite{AndersonDuffin69}. Over $\C$, the rank
equality is the invertible-sum specialization of Tian's formula
\cite[Theorem~10.8(b), pp.~85--86]{Tian00}, which credits Rao and
Mitra in the following discussion. Lemma~\ref{lem:rank} includes
the direct kernel and image proof over an arbitrary field.
The Banach spectral interpretation of generalized Drazin
invertibility is due to Koliha~\cite[Theorems~3.1 and~4.2]{Koliha96}.

Section~\ref{sec:corners} treats coefficient algebras and Hua
configurations. Theorem~\ref{thm:realization} gives index-three
inputs over any nonzero unital ring with $2$ invertible. Its Banach
specialization has spectrum $\{0\}\cup\tfrac12\sigma(Q)$ for any
prescribed coefficient $Q$, after $3\times3$ amplification.
The Volterra and diagonal examples give a quasinilpotent nonnilpotent
residue and an output without a generalized Drazin inverse.
These statements are not replaced by a fixed-dimensional matrix
classification. The latter does yield Corollary
\ref{cor:hua-output-classification}: a finite matrix is a Hua output
in the stated finite-Drazin configuration class exactly when its
nullity is not one. Corollary~\ref{cor:generic-candidate-failure}
specifies the locus on which the zero-core candidate generically
fails, although the output's Drazin inverse always exists.

The companion~\cite{HuaI} establishes the common-core reduction
and inverse criteria. Its Examples~7.1--7.5 overlap with our
scalar-core families, operator boundary instances and weighted-cycle
pairs; here their required proofs support realization and splitting.
Proofs of the general coefficient construction, ring/Banach
scalar-shift identities and separated-cluster upper bound are included here.
Section~\ref{sec:scope} compares eigenvalue rates, spectral gaps and
inverse behavior. Explicit families admit a constant bounded
reflexive inverse while nonzero eigenvalues tend to zero, and the
scalar coefficient family exhibits discontinuous Drazin inversion.
For complex matrices, the finite-dimensional form of
Badea--Mbekhta's gap formula~\cite[Theorem~3.2]{BadeaMbekhta97}
requires the admissible inner inverses to preserve the output's nonzero
primary subspace; preservation of the common input core does not suffice.
The section also gives the precise additive comparison.
Section~\ref{sec:problems} isolates the remaining restricted questions.
The auxiliary Appendix~\ref{sec:cover} treats raw square roots,
critical points, analytic selection and the cover on the specified
simple nonzero-spectrum base.

\section{Kernel, image, and index restrictions}\label{sec:structure}
We write $S=N+C$, $\beta=S^{-1}$ and $\nu=\Par(N,C)$. In any unital ring in which $S$ is a unit,
\begin{equation}\label{eq:symmetry}
 \nu=N-N\beta N=N\beta C=C\beta N=C-C\beta C.
\end{equation}
For example, $N\beta C=N\beta(S-N)=N-N\beta N$ and $C\beta N=(S-N)\beta N=N-N\beta N$. The other equality follows by exchanging $N,C$. No commutation of the two inputs is assumed.

For matrices, $\ker$ and $\im$ refer to their actions on column vectors, and $E_{ij}$ denotes a matrix unit. For complex matrices, we use the operator norm induced by the Euclidean norm. The nilpotency index of a nilpotent element is its least positive vanishing power; in a nonzero matrix algebra, $\ind(0)=1$. A regular nilpotent over a field has a single Jordan block of size $m$.
We specify the ground field in the algebraic realization statements;
all analytic perturbation results are over $\C$.

\begin{lemma}[Kernel, image, and rank]\label{lem:rank}
Let $\mathbb F$ be any field and let $N,C\in M_m(\mathbb F)$,
with $S=N+C$ invertible. Then
\begin{gather}
 \ker\nu=\ker N\oplus\ker C,\qquad
 \im\nu=\im N\cap\im C,\label{eq:kernel-image}\\
 \rk\nu=\rk N+\rk C-m.\label{eq:rank}
\end{gather}
In particular, nilpotent inputs give $\rk\nu\le m-2$.
\end{lemma}
Over $\C$, the rank identity is the specialization of~\cite[Theorem 10.8(b), pp.~85--86]{Tian00} mentioned above: the two range hypotheses there are automatic for an invertible square sum. We include the direct kernel argument and the image conclusion.
\begin{proof}
Put $U=S^{-1}N$ and $V=S^{-1}C=I-U$. Then $UV=VU$ and $\nu=SUV$. If $UVx=0$, write $x=Ux+Vx$. The first summand lies in $\ker V$ and the second in $\ker U$. Moreover $\ker U\cap\ker V=\{0\}$ because $U+V=I$. Conversely both kernels are contained in $\ker(UV)$. Since $\ker U=\ker N$ and $\ker V=\ker C$, this proves the kernel and rank identities.

The two ordered expressions $\nu=N\beta C=C\beta N$ put $\im\nu$ in $\im N\cap\im C$. Also $\im N+\im C=\mathbb F^m$, because it contains the image of $S$. The dimension of the intersection is therefore $\rk N+\rk C-m=\rk\nu$, proving equality. Each nilpotent input has rank at most $m-1$, which gives the final bound.
\end{proof}

\begin{proposition}[Sharp nilpotent-output index]\label{prop:index}
Let $m\ge2$, and suppose $N,C\in M_m(\C)$ are nilpotent, $N+C$ is invertible, and $\nu=\Par(N,C)$ is nilpotent. Then $\ind(\nu)\le m-1$. This bound is attained for every $m\ge2$ by
\begin{equation}\label{eq:cycle-inputs}
 N_m=\sum_{i=1}^{m-1}E_{i,i+1},\qquad
 C_m=E_{m1}+\sum_{i=1}^{m-2}E_{i,i+1}.
\end{equation}
If $\ind(\nu)=m-1$, both inputs are regular nilpotents of index $m$, and $\ker\nu$ is the direct sum of their two kernel lines.
\end{proposition}
\begin{proof}
By Lemma~\ref{lem:rank}, the two nonzero input kernels have zero intersection and lie in $\ker\nu$. A nilpotent output thus has at least two Jordan blocks. Its largest block has size at most $m-1$.

Both inputs in~\eqref{eq:cycle-inputs} are regular nilpotents. For $N_m$ the chain is $e_m\mapsto e_{m-1}\mapsto\cdots\mapsto e_1\mapsto0$; for $C_m$ it is $e_{m-1}\mapsto\cdots\mapsto e_1\mapsto e_m\mapsto0$. The sum is the invertible weighted cycle
\[
 S_me_1=e_m,\qquad S_me_m=e_{m-1},\qquad
 S_me_j=2e_{j-1}\quad(2\le j\le m-1).
\]
Its inverse acts by
\[
 S_m^{-1}e_i=\tfrac12e_{i+1}\ (1\le i\le m-2),\qquad
 S_m^{-1}e_{m-1}=e_m,\qquad S_m^{-1}e_m=e_1.
\]
Consequently $N_mS_m^{-1}C_m$ sends $e_j$ to $e_{j-1}/2$ for $2\le j\le m-1$ and kills $e_1,e_m$. Thus
\begin{equation}\label{eq:cycle-residue}
 \Par(N_m,C_m)=\frac12\sum_{i=1}^{m-2}E_{i,i+1},
\end{equation}
which has index $m-1$. At $m=2$ the sums with upper limit zero are empty, $S_2$ is the swap, and the residue is zero of index one.

Finally, if the output index is $m-1$, its Jordan type is $[m-1,1]$, including $[1,1]$ at $m=2$. Its kernel has dimension two. The direct sum in~\eqref{eq:kernel-image} therefore consists of two lines, so each input has one Jordan block.
\end{proof}
The index bound follows from the rank bound: a nilpotent matrix of
rank $r$ has largest Jordan block at most $r+1$. Proposition~\ref{prop:index}
also gives sharpness in every dimension and the regular-input equality
case. Theorem~\ref{thm:complete-realization} shows that the rank bound
characterizes all complex-matrix residues in every dimension.
Problem~\ref{prob:regime} concerns the further geometry of realizing pairs.

\begin{proposition}[Square-zero inputs]\label{prop:square-zero}
In any unital ring, if $N^2=0$, then $(N-NGN)^2=0$ for every $G$. In particular, one square-zero input in~\eqref{eq:parallel-map} forces a square-zero residue. If both $N^2=C^2=0$ and $N+C$ is a unit, then $\nu=0$.
\end{proposition}
\begin{proof}
Expanding $(N-NGN)^2$ gives four terms, each containing $N^2$, so the square vanishes. Equation~\eqref{eq:symmetry} treats either input. If both squares vanish, then $SN=CN=CS$. Multiplying $SN=CS$ on the left and right by $\beta$ gives $N\beta=\beta C$. Hence $\nu=N\beta C=\beta C^2=0$.
\end{proof}

\section{Exact realization}\label{sec:realization}

The rank restriction is sufficient over every infinite field.
We first prove this fixed-dimensional classification, keeping its
weighted-cycle special case explicit. Section~\ref{sec:corners} treats the coefficient-ring construction
with its separate index-three and Banach spectral conclusions.

The construction uses two of the zero-primary chains to accommodate the
invertible primary part. First, a weighted cycle realizes the nilpotent
target after its closing edges are divided between the two inputs. We
then split the inverse of the invertible part into two cyclic units and
make their inverses nilpotent by corrections of rank at most one. Joining
these corrections to two chain ends preserves nilpotence of the inputs.
The final ordered factorization verifies the prescribed parallel sum.

\begin{theorem}[Complete realization over an infinite field]
\label{thm:complete-realization}\label{thm:M3}
\label{thm:nilpotent-realization}\label{cor:M4-nilpotent}
\label{prob:types}
Let $\mathbb F$ be an infinite field, $m\ge2$ and $T\in M_m(\mathbb F)$.
There exist nilpotent matrices $N,C\in M_m(\mathbb F)$ such that
$N+C$ is invertible and
\[
 T=N(N+C)^{-1}C
\]
if and only if $\rk T\le m-2$.
In that case one input can always be chosen regular nilpotent. More
generally, every pair of positive integers $p,q$ with
$p+q=\dim\ker T$ occurs as
$(\dim\ker N,\dim\ker C)$ for a realizing nilpotent pair.
In every realizing pair, both inputs are regular nilpotent if and only
if $\dim\ker T=2$.
\end{theorem}
\begin{proof}
The field-independent kernel and rank identities of
Lemma~\ref{lem:rank} give necessity and, for every realizing pair,
\[
 \dim\ker T=\dim\ker N+\dim\ker C.
\]
For sufficiency, the Fitting decomposition
$\mathbb F^m=\ker T^m\oplus\im T^m$ and nilpotent Jordan form over
$\mathbb F$ allow a simultaneous similarity putting the target in the form
\[
 T=J\oplus M,\qquad
 J=J_{n_1}\oplus\cdots\oplus J_{n_j},\qquad
 M\in GL_k(\mathbb F),\qquad j=\dim\ker T\ge2.
\]
Here $d=n_1+\cdots+n_j$, $d+k=m$, and the invertible part is absent
when $k=0$. Let $u_i,t_i\in\mathbb F^d$ be the first and last standard
basis columns of the $i$th nilpotent block. Thus $Ju_i=0$ and $t_i^TJ=0$.
All transposes in this proof are ordinary bilinear transposes. Put
\[
 K=\sum_{i=1}^{j-1}t_i u_{i+1}^T,\qquad L=t_j u_1^T.
\]
Choose $t\in\mathbb F\setminus\{0,-1\}$ and set
\[
 n=1+t,\qquad c=1+t^{-1}.
\]
These scalars satisfy
\[
 n,c,n+c\ne0,\qquad n+c=nc,\qquad n^{-1}+c^{-1}=1.
\]
The choice is available because $\mathbb F$ is infinite; it also works
in characteristic two. When the characteristic is not two, one may
take $t=1$ and recover $n=c=2$.

\emph{The weighted-cycle case $k=0$.}
Set
\begin{equation}\label{eq:weighted-realization}
 N=nJ+L,\qquad C=cJ+K,\qquad S=N+C=(n+c)J+K+L.
\end{equation}
The matrix $C$ is strictly upper triangular. The matrix $N$ follows
the directed cycle with positions $E_{i,i+1}$, $1\le i<m$, and
$E_{m1}$, but omits the $j-1\ge1$ cut edges in $K$.
Every walk of length $m$ therefore meets a missing edge, so $N^m=0$.
In $S$, the internal nilpotent-block edges have weight $n+c$, and
the cut and closing edges have weight $1$. Hence $S$ is invertible,
with $\det S=(-1)^{m-1}(n+c)^{m-j}$.

For a cyclic edge with respective weights $a_i,c_i,s_i$ in
$N,C,S$, the ordered product $NS^{-1}C$ has weight $a_ic_i/s_i$
on that edge. This is $nc/(n+c)=1$ inside each nilpotent block,
$0$ at a cut, and $0$ at the closing edge. Thus $NS^{-1}C=J$.
If all blocks have size one, the pair is
$N=E_{m1}$ and $C=\sum_{i=1}^{m-1}E_{i,i+1}$, with zero residue.
This includes the swap pair when $m=2$.

\emph{Cyclic matrices and a rank-one correction.}
Assume now $k\ge1$. A matrix $A\in M_k(\mathbb F)$ is cyclic if it
has a column $v$ such that $v,Av,\ldots,A^{k-1}v$ is a basis.
The cyclic characteristic-polynomial prescription is a standard
rank-one perturbation mechanism; compare
\cite[Proposition~3]{Merzel2021} for the complex invertible case and
nilpotent target used here. We include the elementary argument over
the ground field.

For any cyclic column $v$, there is a row $\xi$ such that
$A+v\xi$ is nilpotent. Write
$\chi_A(z)=z^k+a_1z^{k-1}+\cdots+a_k$. The coefficient columns
of $\adj(zI-A)v$, in descending powers of $z$, are
\[
 v,\quad (A+a_1I)v,\quad\ldots,\quad
 (A^{k-1}+a_1A^{k-2}+\cdots+a_{k-1}I)v.
\]
For $k=1$ the list consists only of $v$.
They form a triangular change of basis from
$v,Av,\ldots,A^{k-1}v$. Consequently
$\xi\mapsto\xi\adj(zI-A)v$ is an isomorphism from rows to
polynomials of degree at most $k-1$. Choose $\xi$ with image
$\chi_A(z)-z^k$. The rank-one determinant identity gives
\[
 \det(zI-A-v\xi)
 =\chi_A(z)-\xi\adj(zI-A)v=z^k.
\]
This identity follows directly by multilinearity in the columns,
and Cayley--Hamilton proves nilpotence. The correction may be zero
when $A$ is already nilpotent; it is nonzero for the invertible
matrices to which we apply it.

We next choose invertible cyclic $B_1,B_2\in M_k(\mathbb F)$ such
that $B_1+B_2=M^{-1}$. Put $D_0=M^{-1}$, and let $\Delta(B)$
denote the discriminant of the characteristic polynomial of $B$,
with $\Delta(B)=1$ when $k=1$. In the entries of a variable matrix
$B$, the expression
\[
 Q(B)=\det B\,\det(D_0-B)\,\Delta(B)\,\Delta(D_0-B)
\]
is a polynomial over $\mathbb F$. It is not the zero polynomial.
Indeed, over an algebraic closure $\overline{\mathbb F}$, triangularize
$D_0$ with diagonal $\mu_1,\ldots,\mu_k$. In those coordinates
choose $B=\operatorname{diag}(z_1,\ldots,z_k)$ so that the $z_i$
and $\mu_i-z_i$ are nonzero and each of these two lists has distinct
entries. The excluded choices form finitely many proper affine
hyperplanes over the infinite field $\overline{\mathbb F}$, so
such a choice exists. Conjugating back gives a matrix at which $Q$
is nonzero. This argument is only a witness that $Q$ is a nonzero
polynomial; the following choice is made over the ground field.

A nonzero polynomial over an infinite field cannot vanish at every
point of that field, by induction on its number of variables.
Choose $B_2\in M_k(\mathbb F)$ with $Q(B_2)\ne0$, and set
$B_1=D_0-B_2$. Both $B_i$, and hence both $A_i=B_i^{-1}$,
have simple nonzero spectrum over $\overline{\mathbb F}$.
They have cyclic columns over $\mathbb F$: for such an $A_i$,
\[
 v\longmapsto\det(v,A_iv,\ldots,A_i^{k-1}v)
\]
is a polynomial over $\mathbb F$ which is nonzero over
$\overline{\mathbb F}$ by the Vandermonde determinant in an
eigenbasis. It therefore has a nonzero value at an $\mathbb F$-column.
Choose such columns $c_i$ and rows $\xi_i$ for which
$D_i=A_i+c_i\xi_i$ is nilpotent, for $i=1,2$.
No commutation of $A_1,A_2$ or $M$ is assumed.

\emph{The input pair.}
With respect to $\mathbb F^d\oplus\mathbb F^k$, define
\begin{equation}\label{eq:complete-inputs}
 N=\begin{pmatrix}
  nJ+L&t_j\xi_1\\ c_1u_1^T&A_1+c_1\xi_1
 \end{pmatrix},\qquad
 C=\begin{pmatrix}
  cJ+K&t_1\xi_2\\ c_2u_2^T&A_2+c_2\xi_2
 \end{pmatrix}.
\end{equation}
For $N$, order the large blocks as the $j$th nilpotent block,
$\mathbb F^k$, the first nilpotent block, and the remaining ones.
For $C$, use the first nilpotent block, $\mathbb F^k$, and then
the second through $j$th nilpotent blocks. In each case the matrix
is block upper triangular with nilpotent diagonal blocks $nJ_{n_i}$
or $cJ_{n_i}$ and $D_1$ or $D_2$. Hence both inputs are nilpotent.

\emph{Auxiliary projections and the sum.}
Put $S=N+C$. Define rows and columns on $\mathbb F^d\oplus\mathbb F^k$ by
\[
 \ell_i=(t_i^T,0),\qquad
 f_j=\binom{t_j}{c_1},\qquad
 f_1=\binom{t_1}{c_2},\qquad
 f_i=\binom{t_i}{0}\quad(1<i<j).
\]
Since $j\ge2$, the special indices $1,j$ are distinct and
$\ell_i f_h=\delta_{ih}$. Thus $P_i=f_i\ell_i$ are pairwise
orthogonal idempotents. Set $P_0=I-\sum_{i=1}^jP_i$.
These are auxiliary projections for this construction. Define
\[
 G=nI_d\oplus A_1M^{-1},\qquad
 G'=cI_d\oplus A_2M^{-1}.
\]
Their inverses are
\[
 G^{-1}=n^{-1}I_d\oplus MB_1,\qquad
 G'^{-1}=c^{-1}I_d\oplus MB_2,
 \qquad G^{-1}+G'^{-1}=I.
\]
The inverses commute because one is $I$ minus the other. Hence
$G,G'$ commute and $GG'=G+G'$. The terminal rows of $S$ are
\begin{align*}
 \ell_jS&=(u_1^T,\xi_1),&
 \ell_1S&=(u_2^T,\xi_2),\\
 \ell_iS&=(u_{i+1}^T,0)\quad(1<i<j).
\end{align*}
Subtracting $\sum_i f_i\ell_iS$ from $S$ gives
\begin{equation}\label{eq:realization-projector-sum}
 P_0S=\operatorname{diag}((n+c)J,A_1+A_2)=GG'T.
\end{equation}
The top block uses $nc=n+c$, while the bottom block uses the
ordered calculation
\[
 A_1M^{-1}A_2=A_1(B_1+B_2)A_2=A_2+A_1.
\]

If $Sx=0$, equation~\eqref{eq:realization-projector-sum} and
invertibility of $G,G'$ give $Tx=0$. Thus
$x=(\sum_i a_i u_i,0)$. The terminal rows give
$a_1=0$ and $a_{i+1}=0$ for $i<j$. Hence $x=0$, so $S$ is
invertible. This establishes the needed inverse before using it.

\emph{The ordered residue calculation.}
Put $D=G'^{-1}$ and $E=P_j+DP_0$. Direct subtraction and
equation~\eqref{eq:realization-projector-sum} give
\[
 N-P_jS=\operatorname{diag}(nJ,A_1)=GT=DP_0S.
\]
Therefore $N=ES$ and $C=(I-E)S$. For every $i$,
$\ell_iD=c^{-1}\ell_i$, so
\[
 P_jDP_0=0,\qquad P_0P_j=0,\qquad P_0DP_0=DP_0.
\]
Consequently $E^2=P_j+D^2P_0$, and
\[
 E(I-E)=D(I-D)P_0=DG^{-1}P_0.
\]
No commutation of $D$ with $P_0$ has been used. It follows that
\[
 NS^{-1}C=E(I-E)S
 =DG^{-1}P_0S=DG^{-1}GG'T=T.
\]
Simultaneous conjugation of the two inputs realizes the original target.

\emph{Input Jordan types and prescribed kernel dimensions.}
In the constructed coordinates,
\begin{equation}\label{eq:realization-input-kernels}
 \ker C=\langle u_1\rangle,\qquad
 \ker N=\langle u_2,\ldots,u_j\rangle,
\end{equation}
where the columns are embedded in $\mathbb F^d\oplus\mathbb F^k$
when $k>0$. For example, $C(x,y)=0$ gives
$u_2^Tx+\xi_2y=0$ and $u_{i+1}^Tx=0$ for $1<i<j$ from its
terminal rows. Its bottom equation then reduces to $A_2y=0$,
so $y=0$; the top equation gives $Jx=0$. This leaves precisely
$\langle u_1\rangle$. The same argument with the terminal row
$(u_1^T,\xi_1)$ proves the assertion for $N$.
For $k=0$ the kernels follow directly from~\eqref{eq:weighted-realization}.
Thus $C$ is regular nilpotent. The restriction of $N$ to the first
and last zero-primary blocks together with the invertible-part
coordinates is nilpotent with a one-dimensional kernel, and hence
is a single Jordan block. The middle blocks are unchanged, so the
block lengths of $N$ are
\[
 n_1+n_j+k,\ n_2,\ldots,n_{j-1}.
\]
The identity $\dim\ker T=\dim\ker N+\dim\ker C$ shows that
both inputs of any realizing pair are regular if and only if
$\dim\ker T=2$.

Finally, let $p,q\ge1$ with $p+q=j$. The construction above gives
the split $(j-1,1)$; exchanging the inputs gives $(1,j-1)$ because
the parallel sum is symmetric. If $p,q\ge2$, divide the $j$ zero
Jordan blocks into groups of $p$ and $q$ blocks and attach the
invertible primary part to either group. Apply the construction
to these two target summands, using kernel splits $(p-1,1)$ and
$(1,q-1)$, respectively. Their direct sum realizes the original
target, after a block permutation, and has input kernel dimensions
$p$ and $q$. Both sums remain invertible. This proves every
admissible positive split.
\end{proof}

\begin{remark}[A finite-field boundary]\label{rem:realization-F2}
Infiniteness cannot simply be omitted from the theorem. Over
$\mathbb F_2$, no nonzero nilpotent rank-one matrix in $M_3(\mathbb F_2)$
is a parallel sum of nilpotent inputs with invertible sum.
The incidence and normal-form reduction used here works over any field
$\mathbb F$. Let $N,C\in M_3(\mathbb F)$ be nilpotent with invertible
sum and nonzero nilpotent residue $\nu=\Par(N,C)$. By
Lemma~\ref{lem:rank}, both inputs are regular nilpotent and
$\rk\nu=1$, so $\nu^2=0$.
Let $K_N,K_C$ be their kernel lines and $W_N,W_C$ their image planes.
Each kernel line lies in its own image plane, and
\[
 \ker\nu=K_N\oplus K_C,\qquad
 \im\nu=W_N\cap W_C\subseteq\ker\nu.
\]
The planes are distinct because $N+C$ is invertible. Exactly one
cross-incidence $K_N\subseteq W_C$ or $K_C\subseteq W_N$ holds.
Indeed, if both held the planes would both equal $K_N\oplus K_C$.
If neither held, a vector in their intersection written as the sum
of a vector in $K_N$ and one in $K_C$ would have both components
zero, contradicting $\rk\nu=1$.

After exchanging the inputs, assume $K_N\subseteq W_C$. Then
$W_C=K_N\oplus K_C$ and $K_C\nsubseteq W_N$. Thus
$W_N\cap W_C=K_N$, so $\im\nu=K_N$; the exchanged case gives
$\im\nu=K_C$. Put $N=J_3$.
A generator of $K_C$ has nonzero third coordinate, so an invertible
polynomial in $J_3$ sends $e_3$ to that generator. Changing basis
by this polynomial keeps $N=J_3$ and makes $K_C=\langle e_3\rangle$;
it also gives $W_C=\langle e_1,e_3\rangle$. The zero third column,
zero second row and zero trace of $C$ therefore give
\[
 C=uE_{12}+vE_{31}+wE_{32}.
\]
Rank two forces $uv\ne0$, and invertibility of $N+C$ forces
$v(1+u)\ne0$. Over $\mathbb F_2$ this is impossible, since
$u\ne0$ implies $u=1$. Thus, for example, the rank-one target
$E_{12}$ is not realizable over $\mathbb F_2$.
This obstruction does not classify the other finite fields.
\end{remark}

We return to complex matrices in the examples that follow.

\begin{example}[A nonnilpotent target]\label{ex:nonnilpotent-realization}
The matrices
\[
 N=\frac16\begin{pmatrix}
 0&0&0&0\\6&0&-1&4\\6&0&2&4\\6&0&-1&-2
 \end{pmatrix},\qquad
 C=\frac16\begin{pmatrix}
 0&6&4&-1\\0&0&0&0\\0&6&-2&-1\\0&6&4&2
 \end{pmatrix}
\]
satisfy $N^4=C^4=0$, $\det(N+C)=-1/4$, and
\[
 \Par(N,C)=\operatorname{diag}(0,0,1,1).
\]
These identities follow by direct multiplication. The target has
rank two and an invertible primary part with a repeated eigenvalue.
\end{example}

Direct sums preserve realization. If $N_i,C_i$ are nilpotent and
$N_i+C_i$ is invertible for $i=1,2$, then
\begin{equation}\label{eq:direct-sum}
 \Par(N_1\oplus N_2,C_1\oplus C_2)
 =\Par(N_1,C_1)\oplus\Par(N_2,C_2),
\end{equation}
because the inverse of the direct sum is the direct sum of the
inverses. The realizable nilpotent targets in dimension four have types
$[3,1]$, $[2,2]$, $[2,1,1]$ and $[1,1,1,1]$.
Proposition~\ref{prop:index} and Example~\ref{ex:ranktwo} give
concrete realizations of the first two. Another realization of
$[2,1,1]$ is
\[
 N=E_{12}+E_{23}+E_{34},\qquad C=E_{41}+E_{24}.
\]
Here $N$ has index four, while $C^2=E_{21}\ne0$ and $C^3=0$.
Their sum has determinant $-1$, and
\[
 (N+C)^{-1}=\begin{pmatrix}
 0&0&0&1\\1&0&0&0\\0&1&-1&0\\0&0&1&0
 \end{pmatrix},\qquad \Par(N,C)=E_{24}.
\]
The zero residue in dimension four is also obtained from two
copies of the pair $E_{12},E_{21}\in M_2(\C)$, whose sum is
invertible and whose residue is zero by
Proposition~\ref{prop:square-zero}.

\section{Differential structure and characteristic coefficients}\label{sec:mechanism}
Throughout the matrix perturbation discussion, paths are analytic in a
complex parameter near zero. A \emph{paired-input path} consists of
$A_\varepsilon,C_\varepsilon$ with invertible sum. Its equivalent
\emph{raw coordinates} are
\[
 a_\varepsilon=A_\varepsilon,\qquad
 b_\varepsilon=(A_\varepsilon+C_\varepsilon)^{-1},\qquad
 T_\varepsilon=a_\varepsilon-a_\varepsilon b_\varepsilon a_\varepsilon.
\]
At a base $(A_0,C_0)=(N,C)$, write $S=N+C$, $\beta=S^{-1}$ and
$\nu=N\beta C$. If $X=A'_0$, $Z=C'_0$ and $Y=b'_0$, then
\begin{equation}\label{eq:path-dictionary}
 Z=-SYS-X,\qquad Y=-\beta(X+Z)\beta.
\end{equation}
Thus arbitrary raw first derivatives $(X,Y)$ are equivalent to arbitrary
paired first derivatives $(X,Z)$. Nilpotence assumptions below concern
the base inputs; nearby inputs need not remain nilpotent.

The following restrictions and choices of extra data can overlap.
\begin{center}
\begin{tabular}{@{}p{0.28\linewidth}p{0.66\linewidth}@{}}
\hline
Restriction or extra data & Requirement\\ \hline
Fixed coupling & $b_\varepsilon=\beta$, equivalently
$A_\varepsilon+C_\varepsilon=S$; the inputs need not be affine.\\[3pt]
Affine fixed coupling & $A_\varepsilon=N+\varepsilon X$,
$C_\varepsilon=C-\varepsilon X$; the base and $X$ determine the whole output path.\\[3pt]
Nilpotent preservation & Both inputs remain nilpotent along the path.\\[3pt]
Augmented outer data & Selected commuting outer inverses
$a'_\varepsilon,c'_\varepsilon$ are added, as in Section~\ref{sec:corners}.\\[3pt]
Finite-Drazin Hua lift & The augmented data also satisfy the nilpotent
input-remainder and intertwining conditions of Section~\ref{sec:corners}.\\ \hline
\end{tabular}
\end{center}
A raw path need not admit a continuous finite-Drazin Hua lift through
prescribed limiting inverse data. The commuting outer-inverse equations
alone are weaker; Example~\ref{ex:twothirds} gives an analytic outer-data
lift but no continuous finite-Drazin lift through the same limiting
inverse data. We use $w^T$ for an ordinary
transpose, so all vector pairings below are bilinear, not Hermitian.

\begin{proposition}[The paired derivative and its image]\label{prop:paired-derivative}
Let $N,C\in M_m(\C)$, assume $S=N+C$ is invertible, and put
$\beta=S^{-1}$ and $\nu=N\beta C$. Nilpotence is not assumed.
For arbitrary paired directions $X,Z$, the output derivative is
\begin{equation}\label{eq:paired-derivative}
 L(X,Z):=D\Par_{(N,C)}[X,Z]
 =C\beta X\beta C+N\beta Z\beta N.
\end{equation}
Write $K_N=\ker N$, $K_C=\ker C$, and let
$L_N=\{\ell:\ell N=0\}$ and $L_C=\{\ell:\ell C=0\}$
be the left nullspaces, consisting of row vectors. With respect
to the bilinear trace pairing $\langle G,H\rangle=\tr(GH)$,
\begin{equation}\label{eq:paired-annihilator}
 (\im L)^\perp
 =\{G:NGN=0=CGC\}
 =(K_N\otimes L_C)\oplus(K_C\otimes L_N).
\end{equation}
Here $k\otimes\ell$ denotes the matrix $k\ell$. Equivalently,
\begin{equation}\label{eq:paired-image}
 \im L=\{H:\ell_C Hk_N=0,\ \ell_N Hk_C=0
 \text{ for all }k_N,k_C,\ell_N,\ell_C\},
\end{equation}
where the four vectors belong to their indicated nullspaces.
If $p=\dim K_N$ and $q=\dim K_C$, then
\begin{equation}\label{eq:paired-corank}
 \operatorname{codim}\im L=2pq.
\end{equation}
For $H=L(X,Z)$ the two remaining kernel pairings satisfy
\begin{equation}\label{eq:paired-self-pairings}
 \ell_N Hk_N=\ell_N Xk_N,\qquad
 \ell_C Hk_C=\ell_C Zk_C.
\end{equation}
Finally, a nonzero rank-one matrix $G\in(\im L)^\perp$
satisfies exactly one of
\[
 NG=0=GC,\qquad CG=0=GN.
\]
\end{proposition}
\begin{proof}
The product rule gives
\[
 L(X,Z)=X\beta C+N\beta Z-N\beta(X+Z)\beta C.
\]
Using $I-N\beta=C\beta$ and $I-\beta C=\beta N$
gives~\eqref{eq:paired-derivative} in the displayed order.
Trace cyclicity then gives
\[
 \tr(GL(X,Z))
 =\tr(X\beta CGC\beta)+\tr(Z\beta NGN\beta).
\]
The trace pairing is nondegenerate, so its vanishing for every
$X,Z$ is equivalent to $CGC=NGN=0$.

Let $I_0=\im N\cap\im C=\im\nu$, using Lemma~\ref{lem:rank}.
Since $S$ is invertible, $\im N+\im C=\C^m$ and
$K_N\cap K_C=0$. A matrix in the annihilator maps
$\im N$ into $K_N$ and $\im C$ into $K_C$, and hence vanishes
on $I_0$. Conversely, such maps are free on complements of
$I_0$ in $\im N$ and $\im C$. Those complements have dimensions
$q$ and $p$, respectively, by~\eqref{eq:rank}. The annihilator
therefore has dimension $pq+qp$.

Each of the two tensor spaces in~\eqref{eq:paired-annihilator}
lies in the annihilator. Their intersection is zero because their
images lie in the disjoint spaces $K_N,K_C$, and their dimensions
are both $pq$. This proves~\eqref{eq:paired-annihilator};
nondegeneracy of the trace pairing gives~\eqref{eq:paired-image}
and~\eqref{eq:paired-corank}.
For~\eqref{eq:paired-self-pairings}, use
$\ell_NC\beta=\ell_N$, $\beta Ck_N=k_N$ and their exchanged
identities in~\eqref{eq:paired-derivative}.

The row spaces $L_N,L_C$ also intersect only at zero, since a
row annihilating $N$ and $C$ annihilates $S$. Relative to bases
adapted to these independent row and column spaces, the two
summands in~\eqref{eq:paired-annihilator} occupy disjoint diagonal
blocks. Their ranks add. A rank-one element therefore belongs to
exactly one summand, giving the ordered alternatives.
\end{proof}

The image in~\eqref{eq:paired-image} concerns unrestricted paired
first derivatives, equivalently unrestricted raw derivatives
by~\eqref{eq:path-dictionary}. Restricting to fixed coupling gives
$Z=-X$; requiring nearby inputs to remain nilpotent imposes further
conditions. No choice of augmented inverse data is made here.
The next identity supplies additional information of order two
along analytic paths through nilpotent inputs.

\begin{proposition}[Adjugates and the last two characteristic coefficients]\label{prop:adjugate}
Let $A,C\in M_m(\R)$, where $\R$ is a commutative unital ring,
$m\ge2$, and $S=A+C$ is invertible. For $T=AS^{-1}C$,
\begin{equation}\label{eq:adjugate-parallel}
 \adj T=\frac{\det A\,\adj C+\det C\,\adj A}{\det S}.
\end{equation}
Now let $A_\varepsilon,C_\varepsilon\in M_m(\C)$ be analytic paths
through nilpotent matrices $N,C$ with invertible sum, and put
$T_\varepsilon=A_\varepsilon(A_\varepsilon+C_\varepsilon)^{-1}
C_\varepsilon$. Its last two characteristic coefficients, up to
their signs, satisfy
\begin{equation}\label{eq:last-two-coefficients}
 \tr(\adj T_\varepsilon)=O(\varepsilon^2),\qquad
 \det T_\varepsilon=O(\varepsilon^2).
\end{equation}
Write $\chi_\nu(z)=z^2p_\nu(z)$ for $\nu=T_0$; this factorization
exists by Lemma~\ref{lem:rank}. For $H=T'_0$,
\begin{equation}\label{eq:polynomial-cancellation}
 \tr\bigl(p_\nu(\nu)H\bigr)=0.
\end{equation}
In particular, if $\nu$ is nilpotent, then
\begin{equation}\label{eq:general-cancellation}
 \tr(\nu^{m-2}H)=0.
\end{equation}
These identities apply to every raw direction: with $\beta=(N+C)^{-1}$,
\begin{equation}\label{eq:raw-derivative-general}
 H=H(X,Y):=X-X\beta N-N\beta X-NYN,
 \qquad X,Y\in M_m(\C).
\end{equation}
\end{proposition}
\begin{proof}
Adjugation reverses products and
$\adj(S^{-1})=(\det S)^{-1}S$. Hence
\[
 \adj T=\frac{1}{\det S}\adj C\,(A+C)\,\adj A
 =\frac{\det A\,\adj C+\det C\,\adj A}{\det S},
\]
using $A\adj A=(\det A)I$ and $\adj C\,C=(\det C)I$
in the displayed order. Only $S$ is inverted.

For an $m\times m$ matrix $M$, the coefficient of $\lambda$ in
$\det(\lambda I-M)$ is $(-1)^{m-1}\tr(\adj M)$. Since $N,C$
are nilpotent and $m\ge2$, their determinants and adjugate traces
vanish. Each term in the trace of the numerator
of~\eqref{eq:adjugate-parallel} is therefore a product of two analytic
functions vanishing at zero. Its denominator is nonzero there. Also
\begin{equation}\label{eq:determinant-parallel-path}
 \det T_\varepsilon
 =\frac{\det A_\varepsilon\det C_\varepsilon}
        {\det(A_\varepsilon+C_\varepsilon)}.
\end{equation}
This proves~\eqref{eq:last-two-coefficients}.

Comparing coefficients in
$(zI-\nu)\adj(zI-\nu)=\chi_\nu(z)I$ shows that the coefficient
of $z$ in $\adj(zI-\nu)$ is $p_\nu(\nu)$. The derivative of
$\chi_{T_\varepsilon}(z)$ at zero parameter is
$-\tr(\adj(zI-\nu)H)$. Its coefficient of $z$ vanishes by
\eqref{eq:last-two-coefficients}, proving
\eqref{eq:polynomial-cancellation}. If $\nu$ is nilpotent,
$p_\nu(z)=z^{m-2}$, which gives~\eqref{eq:general-cancellation}.
Finally $C_\varepsilon=b_\varepsilon^{-1}-a_\varepsilon$ turns
the parallel sum into the raw Hua output and gives the derivative
\eqref{eq:raw-derivative-general}. Every $X,Y$ occurs along
$a_\varepsilon=N+\varepsilon X$,
$b_\varepsilon=\beta+\varepsilon Y$, with the latter invertible
near zero.
\end{proof}

\begin{proposition}[The informative zero-primary types]\label{prop:zero-primary}
Let $N,C\in M_m(\C)$ be nilpotent, $m\ge2$, with $N+C$
invertible, and put $\nu=N(N+C)^{-1}C$. Write
$\chi_\nu(z)=z^h r(z)$, $r(0)\ne0$, and
$p_\nu(z)=z^{h-2}r(z)$. Then $p_\nu(\nu)\ne0$ exactly when
the zero-primary Jordan partition of $\nu$ is $[h-1,1]$.
When this partition occurs and $h\ge3$, $p_\nu(\nu)$ has rank one
and exactly one of
\begin{equation}\label{eq:polynomial-annihilation}
 Np_\nu(\nu)=0=p_\nu(\nu)C,
 \qquad
 Cp_\nu(\nu)=0=p_\nu(\nu)N
\end{equation}
holds. If $h=2$, it is $r(0)$ times the rank-two zero-spectral
projection; the conclusion asserted in this case is the trace identity
\eqref{eq:polynomial-cancellation}.
\end{proposition}
\begin{proof}
There are at least two zero Jordan blocks by Lemma~\ref{lem:rank}.
On the nonzero primary part, $r(\nu)=0$ by Cayley--Hamilton.
On the zero primary part, $p_\nu(\nu)=\nu_0^{h-2}r(\nu_0)$,
where $r(\nu_0)$ is invertible. This is nonzero precisely when a
zero Jordan block has size at least $h-1$. With at least two blocks
and total size $h$, this is exactly $[h-1,1]$. For $h=2$ the
asserted projection formula follows directly.

Now assume the partition is $[h-1,1]$ with $h\ge3$, and put
$R=p_\nu(\nu)$. This matrix has rank one.
Equation~\eqref{eq:polynomial-cancellation}, for every paired
direction by~\eqref{eq:path-dictionary}, puts $R$ in the trace
annihilator of the paired derivative. The last assertion of
Proposition~\ref{prop:paired-derivative} gives exactly the two
exclusive alternatives in~\eqref{eq:polynomial-annihilation}.
\end{proof}

For nilpotent $\nu$ and $m\ge3$, the scalar identity
$\tr(\nu^{m-2}H)=0$ is nonvacuous exactly at type $[m-1,1]$.
At $m=2$, $\nu=0$ but $p_\nu(\nu)=I$, so it says $\tr H=0$.
At type $[2,2]$ in dimension four, $\nu^2=0$ makes this scalar
identity vacuous; the full coefficient is treated in
Proposition~\ref{prop:full-cancellation}.

\begin{theorem}[The three-dimensional rank-one mechanism]\label{thm:mechanism}
Let $N,C\in M_3(\C)$ be nilpotent, $S=N+C$ invertible,
$\beta=S^{-1}$, and $\nu=N\beta C\ne0$ nilpotent.
Then $\rk\nu=1$, $\nu^2=0$, and its Jordan type is $[2,1]$:
zero has algebraic multiplicity three and geometric multiplicity two.
Moreover:
\begin{enumerate}
\item Both inputs are regular nilpotents and
$\ker\nu=\ker N\oplus\ker C=\langle u\rangle\oplus\langle v\rangle$.
\item Exactly one of $\ker N\subseteq\im C$ and
$\ker C\subseteq\im N$ holds.
\item In the first case $N\nu=0=\nu C$ and $\im\nu=\ker N$;
in the second case $C\nu=0=\nu N$ and $\im\nu=\ker C$.
The opposite left product is nonzero in each case.
\item If $\nu=v_0w^T$, $P_\ell=N\beta$ and $P_r=\beta N$,
then, for every $X,Y\in M_3(\C)$,
\[
 w^T\bigl(X(I-P_r)-P_\ell X-NYN\bigr)v_0=0,
 \qquad\text{equivalently }\tr(\nu H(X,Y))=0.
\]
\end{enumerate}
\end{theorem}
\begin{proof}
Lemma~\ref{lem:rank} gives rank one and two independent input kernel
lines; thus each input is regular. A nonzero nilpotent matrix of rank
one has the stated square and Jordan type.
The field-independent incidence argument in
Remark~\ref{rem:realization-F2} gives (2) and the corresponding
image-line identities. In the first case,
$\im\nu=\langle u\rangle=\ker N$, so $N\nu=0$.
Proposition~\ref{prop:zero-primary}, here with $p_\nu(\nu)=\nu$,
then gives $\nu C=0$: its other alternative would also give
$C\nu=0$, impossible since $S\nu\ne0$. In particular
$C\nu\ne0$. Exchanging the inputs proves the second case.
The last assertion is~\eqref{eq:general-cancellation} and
$\tr(v_0w^TH)=w^THv_0$.
\end{proof}

\begin{lemma}[A simultaneous normal form]\label{lem:normal-form}
For nilpotent $N,C\in M_3(\C)$ with invertible sum and nonzero
nilpotent residue, an exchange of the inputs if necessary, followed
by simultaneous similarity, gives
\begin{equation}\label{eq:normal-form}
 N=E_{12}+E_{23},\qquad
 C=uE_{12}+vE_{31}+wE_{32},\qquad uv(1+u)\ne0.
\end{equation}
In these coordinates,
\[
 \beta=\begin{pmatrix}
 -w/[v(1+u)]&0&1/v\\
 1/(1+u)&0&0\\0&1&0
 \end{pmatrix},\qquad
 \nu=\frac{u}{1+u}E_{12}.
\]
The fixed-coupling determinant coefficient and the canonical
parameters used below are
\begin{align}\label{eq:normal-form-parameters}
 \kappa(X)&=-\frac{u}{1+u}X_{31}X_{23},&
 \alpha&=\frac{u-1}{u+1},&
 q&=-\frac{w}{v(1+u)},&d&=\frac1{v(1+u)}.
\end{align}
Here $\alpha=3-2\tr(N\beta)$, $q=\tr\beta$, $d=\det\beta$.
\end{lemma}
\begin{proof}
The field-independent normal-form reduction in
Remark~\ref{rem:realization-F2} gives~\eqref{eq:normal-form}
over $\C$, including $uv(1+u)\ne0$.
Direct multiplication gives $\beta$ and $\nu$, and
$\adj N=E_{13}$, $\adj C=uvE_{32}$. Substitution in
\eqref{eq:fixed-kappa} below gives $\kappa$, while the traces and
determinant give the remaining parameters.
\end{proof}

Exchanging $N$ and $C$ replaces $\alpha$ by $-\alpha$ and leaves
$q$ and $d$ unchanged, since $\tr(C\beta)=3-\tr(N\beta)$.
The polynomial $\kappa(X)$ is unchanged by~\eqref{eq:fixed-kappa}.
Exchanging an actual fixed-coupling path also replaces $X$ by $-X$,
which leaves $\kappa$ unchanged because it is quadratic.

In these coordinates, Proposition~\ref{prop:paired-derivative}
describes the full paired derivative image by the two independent
conditions
\[
 H_{21}=H_{33}=0.
\]
Indeed $K_N=\langle e_1\rangle$, $K_C=\langle e_3\rangle$,
$L_N=\langle e_3^T\rangle$ and $L_C=\langle e_2^T\rangle$.
For paired first derivatives $X,Z$, one also has
\[
 H_{31}=X_{31},\qquad H_{23}=Z_{23},\qquad
 \kappa_{\rm raw}=\frac{u}{1+u}H_{23}H_{31},
\]
with $\kappa_{\rm raw}$ as in~\eqref{eq:raw-kappa} below.
Thus both universal linear constraints, and the first determinant
coefficient, are visible in the same normal form.

In particular the two exceptional hyperplanes are distinct.
For a complex linear space $L\subseteq M_3(\C)$,
$\kappa|_L\equiv0$ if and only if
$L\subseteq\{X_{31}=0\}$ or $L\subseteq\{X_{23}=0\}$,
since the polynomial algebra on $L$ has no zero divisors.
A reducible algebraic set may have components in both hyperplanes.
Also $\alpha=0$ exactly when $u=1$, and $\alpha$ is never $1$ or
$-1$. These facts identify the first exceptional locus without
classifying its finer eigenvalue rates.

\section{Matrix eigenvalue splitting}\label{sec:rates}

All paths in this section are analytic in a complex parameter.

For an analytic input path with nilpotent base inputs, invertible input
sum and nilpotent base residue, all nonleading characteristic coefficients
are $O(\varepsilon)$. Proposition~\ref{prop:adjugate} gives the stronger
$O(\varepsilon^2)$ estimate for the coefficients of $\lambda$ and $1$,
from identities along the full path. This information is additional to
the first-order derivative image in Proposition~\ref{prop:paired-derivative}.
The splitting arguments compare these coefficient orders after rescaling
the spectral variable. Theorems~\ref{thm:generic} and~\ref{thm:extremal-splitting}
impose additional nonvanishing conditions for their generic exact laws.

By a \emph{local ramified parameter} we mean $\varepsilon=\zeta^k$
for a positive integer $k$.
Roots of a monic polynomial $P$ with coefficients holomorphic near
$\varepsilon=0$ admit such parametrizations. Indeed, take its
square-free part $Q$ over the field of meromorphic germs. Since
$Q\mid P\mid Q^{\deg P}$, a sufficiently small punctured disk has
no poles in these polynomial identities and no zeros of the
discriminant of $Q$. Thus the roots of $Q$ are precisely the distinct
roots of $P$ and form a finite unbranched cover. A power substitution
kills their permutation monodromy, making each branch single-valued
and holomorphic off zero. The original monic polynomial has bounded
coefficients and hence uniformly bounded roots; removable singularities
extend each branch holomorphically across zero. Retaining the constant
multiplicities on the punctured disk gives a polynomial factorization
that extends across zero. These convergent Puiseux parametrizations
therefore include repeated limiting roots.

For matrices, write
\[
 \gap(T)=d(0,\sigma(T)\setminus\{0\}),\qquad d(0,\varnothing)=+\infty.
\]
Thus the gap excludes permanent zero roots. The estimates below concern
the indicated path; a generic exact law is stated only with its
specified nonvanishing coefficients.

\subsection{The three-dimensional nilpotent boundary}
\begin{theorem}[A two-thirds bound for every analytic raw path]\label{thm:bound}
Let $N,C\in M_3(\C)$ be nilpotent, $N+C$ invertible, and
$\nu=N(N+C)^{-1}C\ne0$ nilpotent. Put $\beta=(N+C)^{-1}$.
For analytic raw paths $a_\varepsilon,b_\varepsilon$ through
$N,\beta$, every eigenvalue of
$T_\varepsilon=a_\varepsilon-a_\varepsilon b_\varepsilon a_\varepsilon$
satisfies
\[
 \lambda_j(\varepsilon)=O(|\varepsilon|^{2/3}).
\]
The constants are for the fixed base and path. In particular, no branch
has a nonzero leading term of order $\varepsilon^{1/2}$.
\end{theorem}
\begin{proof}
Write
\[
 \chi_{T_\varepsilon}(\lambda)
 =\lambda^3+A_1(\varepsilon)\lambda^2
  +A_2(\varepsilon)\lambda+A_3(\varepsilon).
\]
The limiting trace is zero, so $A_1=O(\varepsilon)$.
Proposition~\ref{prop:adjugate} gives $A_2,A_3=O(\varepsilon^2)$.
For $\varepsilon\ne0$, substitute
$\lambda=|\varepsilon|^{2/3}z$ and divide by $|\varepsilon|^2$.
The new coefficients have orders
$O(|\varepsilon|^{1/3})$, $O(|\varepsilon|^{2/3})$ and $O(1)$.
A monic polynomial whose coefficients have modulus at most $B$
has every root in $|z|\le1+B$, by comparison with its leading term
outside that disk. The roots $z$ are therefore bounded, proving the
estimate and excluding a nonzero square-root leading term.
\end{proof}

\begin{theorem}[First-jet dichotomy and generic fixed-coupling exponent]\label{thm:generic}
Let $N,C\in M_3(\C)$ be nilpotent with invertible sum, put
$\beta=(N+C)^{-1}$, and assume $\nu=N\beta C\ne0$ is nilpotent. For analytic paired-input
paths through $N,C$, set $X=A'_0$, $Z=C'_0$ and
 $T_\varepsilon=A_\varepsilon(A_\varepsilon+C_\varepsilon)^{-1}
 C_\varepsilon$. Put
\begin{equation}\label{eq:raw-kappa}
 \kappa_{\rm raw}=\det(\beta)
       \tr(\adj(N)X)\tr(\adj(C)Z).
\end{equation}
If $\kappa_{\rm raw}\ne0$, both inputs are invertible for small
nonzero parameter, and on a local Puiseux parameter the output roots are
\[
 \lambda_j(\varepsilon)
 =\omega_j\kappa_{\rm raw}^{1/3}\varepsilon^{2/3}
  +O(\varepsilon),\qquad \omega_j^3=1.
\]
If $\kappa_{\rm raw}=0$, every output eigenvalue is $O(|\varepsilon|)$.
Moreover, for $H=T'_0$ and $\nu=v_0w^T$,
\begin{equation}\label{eq:kappa-output-derivative}
 \kappa_{\rm raw}=w^T\adj(H)v_0.
\end{equation}
For affine fixed coupling, $A_\varepsilon=N+\varepsilon X$,
$C_\varepsilon=C-\varepsilon X$, the coefficient is
\begin{equation}\label{eq:fixed-kappa}
 \kappa(X)=-\det(\beta)
          \tr(\adj(N)X)\tr(\adj(C)X).
\end{equation}
Its nonzero set is nonempty and Zariski open in the full direction
space $M_3(\C)$: it is the complement of the two hyperplanes
defined by the displayed trace factors.
\end{theorem}
\begin{proof}
The exact factorization~\eqref{eq:determinant-parallel-path} gives
\[
 \det T_\varepsilon
 =\kappa_{\rm raw}\varepsilon^2+O(\varepsilon^3).
\]
Together with Proposition~\ref{prop:adjugate}, this yields
\begin{equation}\label{eq:three-dimensional-coefficients}
 \chi_{T_\varepsilon}(\lambda)
 =\lambda^3+O(\varepsilon)\lambda^2
  +O(\varepsilon^2)\lambda
  -\kappa_{\rm raw}\varepsilon^2+O(\varepsilon^3).
\end{equation}
If $\kappa_{\rm raw}\ne0$, put $\varepsilon=\zeta^3$ and
$\lambda=\zeta^2w$. Division by $\zeta^6$ gives an analytic
equation whose limiting polynomial is $w^3-\kappa_{\rm raw}$.
Its three roots are simple and nonzero, so the analytic
implicit-function theorem gives the stated expansions.
The two input determinants have nonzero first derivatives,
which also gives their punctured invertibility.
If $\kappa_{\rm raw}=0$, put $\lambda=\varepsilon z$ and
divide~\eqref{eq:three-dimensional-coefficients} by $\varepsilon^3$.
The resulting monic cubic has bounded coefficients; the same root
bound gives $\lambda_j=O(|\varepsilon|)$.

For~\eqref{eq:kappa-output-derivative}, multilinearity gives
\[
 \det(\nu+\varepsilon H+O(\varepsilon^2))
 =\varepsilon^2w^T\adj(H)v_0+O(\varepsilon^3).
\]
Terms containing two columns from the rank-one $\nu$ vanish.
Consequently an output term of order two cannot contribute to the
quadratic determinant coefficient. Comparing coefficients proves
the identity.
Finally fixed coupling gives $Z=-X$. Both limiting inputs have
rank two by Theorem~\ref{thm:mechanism}, so both adjugates are
nonzero and the two trace factors are nonzero linear functionals.
Their product is a nonzero polynomial, proving genericity.
\end{proof}

The formula also applies to nonaffine fixed-coupling paths, with
$Z=-X$, and to arbitrary raw paths via~\eqref{eq:path-dictionary}.
On $\kappa_{\rm raw}=0$, the conclusion is an upper bound;
nonzero linear leading terms are not asserted.

\begin{remark}[The relevant perturbation coefficient]\label{rem:coefficient}
The cancelled scalar $\tr(\nu H)$ belongs to the length-two
Jordan chain in type $[2,1]$. Its vanishing does not isolate two
analytic eigenvalues or imply semisimplicity, linear splitting, or a
square-root series: all three branches can interact. The two-thirds
phenomenon is classical~\cite[pp.~805 and~811]{MBO97}; the statements
above identify the cancellation and determinant factorization forced
by the Hua input map. The generic leading data depend on $H$ by
\eqref{eq:kappa-output-derivative}. On the exceptional locus,
higher output terms can be essential, as
Remark~\ref{rem:census} shows. Our implicit-function arguments use
the explicitly identified simple nonzero rescaled roots.
\end{remark}

\subsection{Other zero-primary types and higher dimensions}
\begin{corollary}[Other base types and a higher-dimensional bound]\label{cor:other-rates}
Let $A_\varepsilon,C_\varepsilon\in M_m(\C)$, $m\ge2$, be analytic
through nilpotent $N,C$, with $N+C$ invertible. Put
$\beta=(N+C)^{-1}$, $X=A'_0$, $Z=C'_0$, and let
$T_\varepsilon=\Par(A_\varepsilon,C_\varepsilon)$, $\nu=T_0$.
Define $\kappa_{\rm raw}$ by the formula~\eqref{eq:raw-kappa}
in this dimension.
\begin{enumerate}
\item Suppose zero has algebraic multiplicity exactly two in $\nu$.
Then $\nu\sim0_2\oplus M$ with $M$ invertible. Put $d_0=\det M$,
with $d_0=1$ if $m=2$. The two roots tending to zero are
$O(|\varepsilon|)$. If $\kappa_{\rm raw}\ne0$, they are analytic
and satisfy
\begin{equation}\label{eq:two-zero-linear}
 \lambda_\pm(\varepsilon)
 =\pm\sqrt{-\kappa_{\rm raw}/d_0}\,\varepsilon
  +O(\varepsilon^2).
\end{equation}
If $\kappa_{\rm raw}=0$, the sharper bound is
\begin{equation}\label{eq:two-zero-exceptional}
 \lambda_\pm(\varepsilon)=O(|\varepsilon|^{3/2}).
\end{equation}
The remaining roots stay separated from zero. In particular, for
$m=3$ and nonnilpotent $\nu$, one has $d_0=\tau=\tr\nu\ne0$,
and the remaining root is $\tau+O(\varepsilon)$.
If $\nu=0$, in any dimension all output roots are
$O(|\varepsilon|)$.
\item If $m\ge4$ and $\nu$ is nilpotent, every output root satisfies
\begin{equation}\label{eq:higher-dimensional-bound}
 \lambda_j(\varepsilon)
 =O\bigl(|\varepsilon|^{1/(m-2)}\bigr).
\end{equation}
\end{enumerate}
All estimates concern the specified fixed path. Neither
\eqref{eq:two-zero-exceptional} nor~\eqref{eq:higher-dimensional-bound}
asserts an exact exponent for every such path.
\end{corollary}
\begin{proof}
For (1), Lemma~\ref{lem:rank} gives $\dim\ker\nu\ge2$,
so the zero primary part is semisimple. Its characteristic
factor $\lambda^2$ is coprime to the nonzero-primary factor.
The linearization of multiplication on the two monic factors
is bijective by coprimality, so the analytic implicit-function
theorem gives
\[
 \chi_{T_\varepsilon}(\lambda)
 =g_\varepsilon(\lambda)
   \bigl(\lambda^2+b_1(\varepsilon)\lambda+b_2(\varepsilon)\bigr),
 \qquad g_0(0)=(-1)^{m-2}d_0\ne0.
\]
Here $g_\varepsilon$ is monic of degree $m-2$ and equals $1$
when $m=2$. The exact determinant factorization
\eqref{eq:determinant-parallel-path} gives
\[
 \det T_\varepsilon
 =\kappa_{\rm raw}\varepsilon^2+O(\varepsilon^3),
 \qquad
 b_2(\varepsilon)=\frac{\kappa_{\rm raw}}{d_0}\varepsilon^2
                   +O(\varepsilon^3).
\]
The full coefficient of $\lambda$ is $O(\varepsilon^2)$ by
Proposition~\ref{prop:adjugate}, and equals
$g_\varepsilon(0)b_1+g'_\varepsilon(0)b_2$.
Thus $b_1=O(\varepsilon^2)$.
When $\kappa_{\rm raw}\ne0$, substitute $\lambda=\varepsilon w$
in the quadratic and divide by $\varepsilon^2$. The resulting
analytic equation has the two simple limiting roots
$\pm\sqrt{-\kappa_{\rm raw}/d_0}$, proving
\eqref{eq:two-zero-linear}. When $\kappa_{\rm raw}=0$,
$b_2=O(\varepsilon^3)$, and the quadratic formula gives
\eqref{eq:two-zero-exceptional}. The nonzero cluster stays separated
by continuity. In dimension three its single root is simple,
giving the stated analytic expansion. If $\nu=0$, the estimate
$T_\varepsilon=O(\varepsilon)$ bounds every eigenvalue directly.

For (2), write
\[
 \chi_{T_\varepsilon}(\lambda)
 =\lambda^m+\sum_{k=1}^{m}a_k(\varepsilon)\lambda^{m-k}.
\]
All $a_k=O(\varepsilon)$, and the last two are
$O(\varepsilon^2)$ by Proposition~\ref{prop:adjugate}.
Set $\lambda=|\varepsilon|^{1/(m-2)}z$ and divide by
$|\varepsilon|^{m/(m-2)}$. For $k\le m-2$ the new coefficient
has order $|\varepsilon|^{1-k/(m-2)}$; the last two have orders
$|\varepsilon|^{(m-3)/(m-2)}$ and
$|\varepsilon|^{(m-4)/(m-2)}$. These exponents are nonnegative,
so the resulting monic polynomial has bounded coefficients and
bounded roots.
\end{proof}

The exceptional estimate~\eqref{eq:two-zero-exceptional} is sharp
already for affine fixed coupling in dimension three. Take
\[
 N=E_{12}+E_{23},\qquad
 C=\begin{pmatrix}-2&1&0\\-4&2&0\\1&0&0\end{pmatrix},\qquad
 X=E_{21}+E_{23}+E_{33}.
\]
Both base inputs are regular nilpotent, $\det(N+C)=2$, and
$\nu=\tfrac12E_{12}+E_{22}$. For
$A_\varepsilon=N+\varepsilon X$, $C_\varepsilon=C-\varepsilon X$,
direct multiplication gives
\[
 \det A_\varepsilon=-\varepsilon^2,\qquad
 \det C_\varepsilon=-\varepsilon(1+\varepsilon),
\]
and
\begin{equation}\label{eq:two-zero-sharp-cubic}
 \chi_{T_\varepsilon}(\lambda)
 =\lambda^3+(3\varepsilon^2-\varepsilon-1)\lambda^2
   +\tfrac12\varepsilon^2\lambda
   -\tfrac12\varepsilon^3(1+\varepsilon).
\end{equation}
Thus $\kappa_{\rm raw}=0$. Substituting
$\varepsilon=\zeta^2$, $\lambda=\zeta^3w$ and dividing by $\zeta^6$
gives
\[
 \zeta^3w^3+(3\zeta^4-\zeta^2-1)w^2
       +\tfrac12\zeta w-\tfrac12(1+\zeta^2)=0.
\]
Its simple limiting roots are $w_0=\pm i/\sqrt2$; implicit
differentiation gives $w'(0)=1/4$. Consequently the small roots are
\[
 \lambda_\pm
 =\pm\frac{i}{\sqrt2}\varepsilon^{3/2}
       +\frac14\varepsilon^2+O(\varepsilon^{5/2})
\]
on the local ramified parameter. This example concerns a nonnilpotent
base residue and is distinct from Theorem~\ref{thm:generic}.

If the nilpotency index of $\nu$ is $k$, then
$T_\varepsilon^k=O(\varepsilon)$ gives the elementary bound
$O(|\varepsilon|^{1/k})$. Thus~\eqref{eq:higher-dimensional-bound}
improves that estimate precisely at type $[m-1,1]$; for $k\le m-2$
the index estimate is at least as strong.

\begin{theorem}[Generic splitting at an extremal nilpotent residue]\label{thm:extremal-splitting}
Let $m\ge4$, let $N,C\in M_m(\C)$ be nilpotent with invertible
sum, and assume $\nu=N(N+C)^{-1}C$ has Jordan type $[m-1,1]$.
For analytic paired-input paths $A_\varepsilon,C_\varepsilon$
through $N,C$, put $\beta=(N+C)^{-1}$, $X=A'_0$, $Z=C'_0$,
$T_\varepsilon=\Par(A_\varepsilon,C_\varepsilon)$ and
\[
 H=T'_0,\qquad c_1=\tr(\nu^{m-3}H),\qquad
 \kappa_{\rm raw}
 =\det(\beta)\tr(\adj(N)X)\tr(\adj(C)Z).
\]
Suppose $c_1\kappa_{\rm raw}\ne0$. Both perturbed inputs and
$T_\varepsilon$ are then invertible for small nonzero parameter.
If $m\ge5$, on a local ramified parameter the roots form two groups:
\begin{align}\label{eq:extremal-leading}
 \lambda_j(\varepsilon)
 &=\alpha_j\varepsilon^{1/(m-2)}(1+o(1)),
 &\alpha_j^{m-2}&=c_1,\quad 1\le j\le m-2,\\
 \lambda_\pm(\varepsilon)
 &=\pm\sqrt{(-1)^m\kappa_{\rm raw}/c_1}\,
       \varepsilon^{1/2}(1+o(1)).\nonumber
\end{align}
Here the $\alpha_j$ are the distinct roots of $z^{m-2}-c_1$.
If $m=4$, the leading multiset is
\begin{equation}\label{eq:extremal-quartic}
 \{z_j\varepsilon^{1/2}(1+o(1)):1\le j\le4\},\qquad
 z_j^4-c_1z_j^2+\kappa_{\rm raw}=0,
\end{equation}
with roots counted with multiplicity.

At every such base, $c_1\kappa_{\rm raw}\ne0$ defines a nonempty
Zariski-open set in the unrestricted paired direction space
$M_m(\C)^2$ and also in the fixed-coupling direction space $Z=-X$.
These genericity assertions concern first derivatives; the leading
laws hold for every analytic path with those derivatives.
\end{theorem}

The two groups reflect the leading balances in
\[
 \lambda^m-c_1\varepsilon\lambda^2
       +(-1)^m\kappa_{\rm raw}\varepsilon^2.
\]
For $m\ge5$, balancing the first two terms predicts
$\lambda^{m-2}\sim c_1\varepsilon$, while balancing the last
two predicts $\lambda^2\sim(-1)^m\kappa_{\rm raw}\varepsilon/c_1$.
At $m=4$ the two scales coincide, and all three terms enter the
quartic balance. The proof below controls the remaining coefficients
and accounts for every root.

\begin{proof}
Write $\chi_{T_\varepsilon}(\lambda)
=\lambda^m+\sum_{k=1}^m a_k(\varepsilon)\lambda^{m-k}$.
Since $\nu$ is nilpotent,
\[
 \adj(\lambda I-\nu)=
 \sum_{j=0}^{m-1}\lambda^{m-1-j}\nu^j.
\]
Differentiating the determinant gives
$a'_k(0)=-\tr(\nu^{k-1}H)$. Thus, with $q=m-2$,
\[
 a_q=-c_1\varepsilon+O(\varepsilon^2),\qquad
 a_{m-1}=O(\varepsilon^2),\qquad
 a_m=(-1)^m\kappa_{\rm raw}\varepsilon^2+O(\varepsilon^3),
\]
where the last two assertions follow from
Proposition~\ref{prop:adjugate} and
\eqref{eq:determinant-parallel-path}. All other coefficients
are $O(\varepsilon)$.

For $m\ge5$, set $\varepsilon=\zeta^{2q}$.
Substitution of $\lambda=\zeta^2w$ and division by $\zeta^{2m}$
give an equation holomorphic at $\zeta=0$ whose limiting
polynomial is
\[
 w^m-c_1w^2=w^2(w^q-c_1).
\]
The $q$ nonzero roots are simple, so the analytic implicit-function
theorem gives $q$ roots with the first asymptotics in
\eqref{eq:extremal-leading}. Next substitute $\lambda=\zeta^qv$
in the original polynomial and divide by $\zeta^{4q}$.
Every term again has a nonnegative power of $\zeta$:
the monic term has exponent $q(q-2)>0$, and a term with
$k<q$ has exponent at least $q(q-k)>0$.
The limiting equation is
\[
 -c_1v^2+(-1)^m\kappa_{\rm raw}=0.
\]
Its two roots are simple and nonzero, giving the remaining pair.
The two scales are distinct because $q>2$, and these roots
account for the degree $m$. For $m=4$, use
$\varepsilon=\zeta^2$, $\lambda=\zeta w$ and divide by $\zeta^4$.
The resulting monic polynomial is holomorphic at zero and tends to
$w^4-c_1w^2+\kappa_{\rm raw}$. Root continuity gives
\eqref{eq:extremal-quartic}. The opening ramification argument
supplies local Puiseux parametrizations even when the limiting roots
repeat; only the displayed leading multiset and relative $o(1)$
remainder are asserted.

It remains to prove genericity at every base. Lemma~\ref{lem:rank}
forces both input ranks to be $m-1$, so their adjugates are
nonzero. Put $G=\nu^{m-3}$ and $R=\nu^{m-2}$.
By Proposition~\ref{prop:zero-primary}, after exchanging the
inputs if necessary, $NR=0=RC$.
Consequently
\[
 RN=RS\ne0,\qquad RN=\nu GN=C\beta NGN,
\]
so $NGN\ne0$. The trace formula in
Proposition~\ref{prop:paired-derivative} therefore shows that
$c_1=\tr(GL(X,Z))$ is a nonzero paired-direction functional.

For fixed coupling, that functional vanishes identically exactly
when $CGC=NGN$. We show that this equality is impossible.
The rank-one space $\im R$ equals $\ker N$; choose a nonzero
$x_1$ in it and $x_2\in\im G$ with $\nu x_2=x_1$.
Then $\im G=\langle x_1,x_2\rangle$, $Nx_1=0$, and
$\im NGN=\langle Nx_2\rangle\ne0$.
Since $\nu GC=RC=0$,
\[
 \im GC\subseteq\im G\cap\ker\nu=\langle x_1\rangle,
 \qquad \im CGC\subseteq\langle Cx_1\rangle.
\]
The alleged equality would give $Nx_2=aCx_1$ for some $a\ne0$.
Applying $C\beta$ and using $C\beta N=\nu$,
$C\beta C=C-\nu$ and $\nu x_1=0$ gives
$x_1=aCx_1$. This contradicts nilpotence of $C$.
Thus $c_1$ is nonzero also on the fixed-coupling direction space.

In each of the two direction spaces, the two determinant-derivative
factors in $\kappa_{\rm raw}$ are nonzero linear functionals.
Their product with $c_1$ is a nonzero polynomial over $\C$,
so its nonzero set is nonempty and Zariski open.
Finally $\kappa_{\rm raw}\ne0$ makes the first derivatives of
both input determinants nonzero. This proves punctured
invertibility of the inputs and, by the determinant factorization,
of the output.
\end{proof}

For $m\ge5$, the maximal root modulus in
\eqref{eq:extremal-leading} has two-sided order
$|\varepsilon|^{1/(m-2)}$, whereas the spectral gap has order
$|\varepsilon|^{1/2}$. At $m=4$ both have order
$|\varepsilon|^{1/2}$; the partition here is $[3,1]$,
distinct from the square-zero type $[2,2]$ in
Proposition~\ref{prop:full-cancellation}.
The dimension-three endpoint is Theorem~\ref{thm:generic},
with its $2/3$ law. A path keeping both inputs nilpotent has
$\kappa_{\rm raw}=0$ and retains at least two zero roots;
the genericity assertion above is in the unrestricted raw
direction spaces.

If both inputs remain nilpotent along a path, Lemma~\ref{lem:rank}
gives at least two permanent zero eigenvalues. In dimension three
the characteristic polynomial is exactly
$\lambda^2(\lambda-\tr T_\varepsilon)$, so the only possibly
moving branch is analytic. In every dimension $m\ge4$,
the bound~\eqref{eq:higher-dimensional-bound} is nevertheless
sharp even while both inputs remain regular nilpotent.

Indeed, put
\[
 J=\sum_{i=1}^{m-2}E_{i,i+1},\qquad
 N=2J+E_{m1},\qquad C=2J+E_{m-1,m},
\]
and set $R=E_{m-1,1}$, $U_\varepsilon=I+\varepsilon R$,
$C_\varepsilon=U_\varepsilon C U_\varepsilon^{-1}$.
The matrices $N,C$ are regular nilpotent: their nonzero edges
each form one chain through all $m$ basis vectors. Since $R^2=0$,
$U_\varepsilon^{-1}=I-\varepsilon R$, so $C_\varepsilon$
remains regular nilpotent for every parameter.
Writing $S_0=N+C=4J+E_{m-1,m}+E_{m1}$, direct ordered
multiplication gives
\[
 C_\varepsilon-C
 =2\varepsilon(E_{m-1,2}-E_{m-2,1}),\qquad
 B_\varepsilon:=S_0^{-1}(C_\varepsilon-C)
 =2\varepsilon E_{m2}-\tfrac12\varepsilon E_{m-1,1}.
\]
Here $B_\varepsilon^2=0$, and hence
\[
 (N+C_\varepsilon)^{-1}
 =(I-B_\varepsilon)S_0^{-1},\qquad
 \det(N+C_\varepsilon)=(-1)^{m-1}4^{m-2}.
\]
Also
$S_0^{-1}N=\operatorname{diag}(1,1/2,\ldots,1/2,0)$ and
$NB_\varepsilon=-\varepsilon E_{m-2,1}$. Substitution in
$\Par(N,C_\varepsilon)=N-N(N+C_\varepsilon)^{-1}N$
therefore gives
\begin{equation}\label{eq:nilpotent-preserving-sharpness}
 \Par(N,C_\varepsilon)=J-\varepsilon E_{m-2,1},\qquad
 \chi_{\Par(N,C_\varepsilon)}(\lambda)
 =\lambda^2(\lambda^{m-2}+\varepsilon).
\end{equation}
For the characteristic polynomial, the first $m-2$ coordinates
form an invariant companion block; the quotient has zero
diagonal blocks of total dimension two.
Thus two roots stay zero and the other $m-2$ have exact modulus
$|\varepsilon|^{1/(m-2)}$. At zero parameter the residue has
type $[m-1,1]$. This sharpness family has
$\kappa_{\rm raw}=0$, consistently with the permanent zero
roots; the generic laws in which every root moves leave the
nilpotent-input class.

\emph{A constant reflexive inverse.}
A reflexive inverse of $T$ is a matrix $G$ satisfying
$TGT=T$ and $GTG=G$; commutation is not required.
For $T_\varepsilon=\Par(N,C_\varepsilon)$ in
\eqref{eq:nilpotent-preserving-sharpness}, put $q=m-2$ and
$G_J=J^T=\sum_{i=1}^{q}E_{i+1,i}$. Direct multiplication gives
\[
 T_\varepsilon G_J=P:=\operatorname{diag}(I_q,0,0),\qquad
 PT_\varepsilon=T_\varepsilon,\qquad G_JP=G_J.
\]
Thus $G_J$ is a reflexive inverse of every $T_\varepsilon$.
Moreover $G_J^q=E_{q+1,1}\ne0$, $G_J^{q+1}=0$ and
$\|G_J\|=1$. The sharp fractional eigenvalue scale
$|\varepsilon|^{1/(m-2)}$ therefore coexists with a constant
nilpotent reflexive inverse of norm one.

\emph{A fixed-sum path with every root moving.}
Keep the base $N,C,J$ used in~\eqref{eq:nilpotent-preserving-sharpness},
and put $q=m-2\ge2$. Set
\[
 X=E_{q+1,m}+E_{m1}+E_{q1},\qquad
 \widetilde A_\varepsilon=N+\varepsilon X,\qquad
 \widetilde C_\varepsilon=C-\varepsilon X.
\]
Their sum is fixed and invertible, with
\[
 \beta=(N+C)^{-1}=\tfrac14J^T+E_{m,q+1}+E_{1m},\qquad
 \det(N+C)=(-1)^{m-1}4^q.
\]
The ordered products
$X\beta X=X$, $X\beta C=E_{q+1,m}$ and
$N\beta X=E_{m1}+\tfrac12E_{q1}$ give
\begin{align*}
 \widetilde T_\varepsilon
 &=\Par(\widetilde A_\varepsilon,\widetilde C_\varepsilon)\\
 &=J-\varepsilon(\varepsilon+\tfrac12)E_{q1}
       +\varepsilon(1-\varepsilon)E_{q+1,m}
       -\varepsilon(1+\varepsilon)E_{m1},\\
 \chi_{\widetilde T_\varepsilon}(\lambda)
 &=\lambda^m+(\varepsilon/2+\varepsilon^2)\lambda^2
       +\varepsilon^2-\varepsilon^4.
\end{align*}
For the characteristic polynomial, put
$u=\varepsilon(\varepsilon+1/2)$ and
$B=\lambda I_q-J_q+uE_{q1}$, where
$J_q=\sum_{i=1}^{q-1}E_{i,i+1}$.
Then $\det B=\lambda^q+u$ and
$e_1^TB^{-1}e_q=(\lambda^q+u)^{-1}$ over the rational-function
field. The Schur complement of this upper-left block in
$\lambda I-\widetilde T_\varepsilon$ gives
$\lambda^2(\lambda^q+u)+\varepsilon^2(1-\varepsilon^2)$.
Polynomial identity extends the formula to every $\lambda$.
Expanding the input determinants gives
\[
 \det\widetilde A_\varepsilon
   =(-1)^{m-1}2^q\varepsilon(1+\varepsilon),\qquad
 \det\widetilde C_\varepsilon
   =(-1)^m2^q\varepsilon(1-\varepsilon).
\]
Thus $c_1=\tr(J^{m-3}\widetilde T'_0)=-1/2$ and
$\kappa_{\rm raw}=(-1)^m$ in Theorem~\ref{thm:extremal-splitting}.
Both inputs and the output are invertible for $0<|\varepsilon|<1$;
the inputs are nilpotent only at the base in this neighborhood.

For $m\ge5$, on a local ramified parameter the theorem gives
\[
 \begin{gathered}
 \lambda_j=\alpha_j\varepsilon^{1/q}(1+o(1)),\quad
 \alpha_j^q=-\tfrac12\quad(1\le j\le q),\\
 \lambda_\pm=\pm i\sqrt2\,\varepsilon^{1/2}(1+o(1)).
 \end{gathered}
\]
Consequently
\[
 \max_{\lambda\in\sigma(\widetilde T_\varepsilon)}|\lambda|
       \sim2^{-1/q}|\varepsilon|^{1/q},\qquad
 \gap(\widetilde T_\varepsilon)\sim\sqrt2\,|\varepsilon|^{1/2}.
\]
At $m=4$, the limiting polynomial is $z^4+z^2/2+1$;
all four roots move on the square-root scale.
Thus~\eqref{eq:nilpotent-preserving-sharpness} retains nilpotent
inputs and two zero roots, while this fixed-sum path illustrates
the generic law with every root moving.

\subsection{Exact paths and full cancellation}
\begin{example}[An exact two-thirds path]\label{ex:twothirds}
Let
\[
N=E_{12}+E_{23},\quad C=E_{12}+E_{31},\quad
\beta=\begin{pmatrix}0&0&1\\1/2&0&0\\0&1&0\end{pmatrix},\quad X=E_{23}+E_{31}.
\]
For $A_\varepsilon=N+\varepsilon X$ and
$C_\varepsilon=C-\varepsilon X$, one has $\nu=E_{12}/2$ and
\[
T_\varepsilon=\begin{pmatrix}
0&1/2&0\\0&0&-\varepsilon-\varepsilon^2\\
\varepsilon-\varepsilon^2&0&0
\end{pmatrix},\qquad
\det(\lambda I-T_\varepsilon)=\lambda^3+\frac{\varepsilon^2-\varepsilon^4}{2}.
\]
Also $\det A_\varepsilon=\varepsilon(1+\varepsilon)$ and
$\det C_\varepsilon=\varepsilon(\varepsilon-1)$.
All three eigenvalues therefore have leading exponent $2/3$.
Adjoining a scalar core as in Section~\ref{sec:corners} gives
$a_\varepsilon=1\oplus(N+\varepsilon X)$ and
$b=\tfrac12\oplus\beta$. The limiting configuration has a nonzero
proper core and nonzero nilpotent residue. With the canonical Drazin
inverse coordinates, the punctured configurations are classical, with core
$I$, since both inputs are invertible. These inverse coordinates do not
extend continuously through zero: a continuous extension of the ordinary
inverses would invert the singular limiting inputs. Thus the raw path has
no continuous lift through those limiting data in the finite-Drazin Hua
class. It does have an analytic lift satisfying the commuting outer-inverse
and intertwining equations: take $a'_\varepsilon=c'_\varepsilon=e=1\oplus0$.
Then $a_\varepsilon e=ea_\varepsilon=c_\varepsilon e=ec_\varepsilon=e$
and $eb=be$. For nonzero parameter the complementary input blocks are
invertible, so the remainders of this constant selection are not nilpotent.
\end{example}

\begin{theorem}[Canonical cubic and exact rate alternatives]\label{thm:canonical}
Let $N,C\in M_3(\C)$ be nilpotent with invertible sum,
$\beta=(N+C)^{-1}$, and nilpotent $\nu=N\beta C$, possibly zero. Put
\[
\alpha=3-2\tr(N\beta),\qquad q=\tr\beta,\qquad d=\det\beta.
\]
For $T_\varepsilon=(N+\varepsilon I)\beta(C-\varepsilon I)$,
\[
\det(\lambda I-T_\varepsilon)
=\lambda^3-(\alpha\varepsilon-q\varepsilon^2)\lambda^2+d\varepsilon^6.
\]
If $\alpha=0$, its eigenvalues are exactly $\varepsilon^2$ times the nonzero roots of $w^3+qw^2+d$, counted with multiplicity. If $\alpha\ne0$, their leading exponent list is $[1,5/2,5/2]$.
\end{theorem}
\begin{proof}
For $\varepsilon\ne0$, the input inverses terminate because $N^3=C^3=0$, and the parallel-sum formula gives
\[
T_\varepsilon^{-1}
=-\varepsilon^{-2}(N+C)+\varepsilon^{-3}(N^2-C^2).
\]
Its trace is zero. The coefficient of $\lambda$ in the characteristic polynomial is
$\det(T_\varepsilon)\tr(T_\varepsilon^{-1})=0$.
Direct expansion gives
\[
 T_\varepsilon=\nu+\varepsilon(C\beta-\beta N)-\varepsilon^2\beta.
\]
Taking its trace gives $\tr T_\varepsilon=\alpha\varepsilon-q\varepsilon^2$,
since $\tr\nu=0$ and $C\beta=I-N\beta$. The product of input determinants gives
$\det T_\varepsilon=-d\varepsilon^6$. These formulas prove the cubic, including at zero by polynomial identity.

For $\alpha=0$, substitution $\lambda=\varepsilon^2w$ gives the fixed cubic. Its constant term is nonzero.
For $\alpha\ne0$, substitute $\lambda=\varepsilon w$ and divide by $\varepsilon^3$. The limiting polynomial is $w^2(w-\alpha)$, with a simple root at $\alpha$, giving one eigenvalue $\alpha\varepsilon+O(\varepsilon^2)$.
For the other two set $\varepsilon=\zeta^2$ and $\lambda=\zeta^5w$. Dividing the equation by $\zeta^{12}$ gives
\[
\zeta^3w^3-(\alpha-q\zeta^2)w^2+d=0.
\]
The two roots $w=\pm\sqrt{d/\alpha}$ are simple at zero. For $G(\zeta,w)=\zeta^3w^3-(\alpha-q\zeta^2)w^2+d$, one has $G_\zeta(0,w_0)=0$ and $G_w(0,w_0)\ne0$. Hence the implicit branch has $w'(0)=0$ and $w(\zeta)=w_0+O(\zeta^2)$. The implicit-function theorem therefore gives
$\lambda=\pm\sqrt{d/\alpha}\,\varepsilon^{5/2}(1+O(\varepsilon))$.
\end{proof}

\begin{remark}[Path dependence does not exclude a generic rate]\label{rem:census}
At the limiting pair in Example \ref{ex:twothirds}, the canonical direction $X=I$ gives the cubic $\lambda^3+\varepsilon^6/2$, hence rates $[2,2,2]$. The direction in that example gives $[2/3,2/3,2/3]$ at the same limit, and both paths have nonzero first derivatives as input paths. Thus actual rates depend on direction, while Theorem \ref{thm:generic} determines the generic rate in its specified direction space.

The canonical exponent is not always two. For $N=E_{12}+E_{23}$ and $C=2E_{12}+E_{31}$, the residue is $2E_{12}/3$, but
\[
\det(\lambda I-T_\varepsilon)=\lambda^3-\frac{\varepsilon}{3}\lambda^2+\frac{\varepsilon^6}{3}.
\]
Its canonical rates are $[1,5/2,5/2]$.

Higher output terms can be essential on the exceptional locus. At the first limiting pair, the canonical derivative is $H=\operatorname{diag}(1/2,-1/2,0)$. The affine approximation has
\[
 \det(\lambda I-(\nu+\varepsilon H))=\lambda^3-\frac{\varepsilon^2}{4}\lambda,
\]
and eigenvalues $0,\varepsilon/2,-\varepsilon/2$. The actual output includes $-\varepsilon^2\beta$ and has the cubic $\lambda^3+\varepsilon^6/2$. Here the quadratic output term changes the rates. By contrast,
when $\kappa_{\rm raw}\ne0$, equation~\eqref{eq:kappa-output-derivative}
already determines the leading two-thirds coefficient from $H$.
\end{remark}

\begin{proposition}[Full cancellation and the rank-two boundary]\label{prop:full-cancellation}\label{prob:higher}
Let $N,C\in M_m(\C)$, with $S=N+C$ invertible and $\beta=S^{-1}$,
and suppose $\nu=N\beta C\ne0$ has rank $r$.
Choose $\nu=VW$, where $V$ has full column rank $r$ and $W$ has
full row rank $r$. Set
\[
 H(X,Y)=X-X\beta N-N\beta X-NYN,\qquad
 \Gamma(X,Y)=WH(X,Y)V.
\]
The identity $\Gamma(X,Y)=0$ for every raw direction $X,Y$ is
equivalent to exactly one of the ordered alternatives
\begin{equation}\label{eq:full-annihilation}
 N\nu=0=\nu C,\qquad\text{or}\qquad C\nu=0=\nu N.
\end{equation}
These conditions force $\nu^2=0$ and $m\ge3r$.
The dimension bound is attained by nilpotent input pairs for every
$r\ge1$; the inequality alone is not sufficient for a given pair.
If $\nu^2=0$ is assumed, the fixed-coupling condition
$\Gamma(X,0)=0$ for every $X$ is equivalent to the same alternatives.
No nilpotence of $N,C$ is required for these assertions.

Now assume $\nu^2=0$, $m=4$ and $r=2$. Universal full cancellation
is impossible. For any analytic raw output path through this base,
put $\Gamma=WT'_0V$. Then
\begin{equation}\label{eq:rank-two-determinant}
 \det T_\varepsilon=(\det\Gamma)\varepsilon^2+O(\varepsilon^3).
\end{equation}
If $\det\Gamma\ne0$ and $\gamma_1,\gamma_2$ are the eigenvalues
of $\Gamma$ with algebraic multiplicities, the leading multiset is
\begin{equation}\label{eq:rank-two-leading}
 \lambda_{j,\pm}(\varepsilon)
 =\pm\sqrt{\gamma_j}\,\varepsilon^{1/2}(1+o(1)),\qquad j=1,2,
\end{equation}
on a local ramified parameter. In particular all four eigenvalue
moduli have two-sided order $|\varepsilon|^{1/2}$.

If the dimension-four inputs $N,C$ are nilpotent, the affine
fixed-coupling directions satisfying
\begin{equation}\label{eq:rank-two-kappa}
 \kappa_4(X):=-\det(\beta)
       \tr(\adj(N)X)\tr(\adj(C)X)\ne0
\end{equation}
form a nonempty Zariski-open set. For every fixed-coupling direction,
$\det\Gamma(X,0)=\kappa_4(X)$. On its nonzero set the law
\eqref{eq:rank-two-leading} holds, and both perturbed inputs are
invertible for small nonzero parameter.
\end{proposition}
\begin{proof}
By~\eqref{eq:path-dictionary}, unrestricted raw directions are
equivalent to unrestricted paired directions. Proposition
\ref{prop:paired-derivative} gives
\[
 WL(X,Z)V=WC\beta X\beta CV+WN\beta Z\beta NV.
\]
For compatible matrices $A,B$, the identity $AXB=0$ for every
matrix $X$ holds exactly when $A=0$ or $B=0$: if both are nonzero,
choose a column $u$ and a row $v$ with $Au\ne0$, $vB\ne0$ and
use $X=uv$. Consequently universal cancellation is equivalent to
\[
 (WC=0\text{ or }CV=0),\qquad
 (WN=0\text{ or }NV=0).
\]
The possibilities $WC=WN=0$ and $CV=NV=0$ would give
$WS=0$ and $SV=0$, respectively, and are excluded by
invertibility and full rank. The remaining two possibilities
are exactly~\eqref{eq:full-annihilation}. They are exclusive,
since both would imply $S\nu=0$. Each forces $\nu^2=0$
by the ordered identities~\eqref{eq:symmetry}.
Full rank of $V,W$ then gives $WV=0$.
The cancellation is independent of the factorization, since
$WHV=0$ is equivalent to $\nu H\nu=0$.

For the fixed-coupling converse, suppose $\nu^2=0$ and
$\Gamma(X,0)=0$ for every $X$. Full rank gives $WV=0$ and
$\nu V=0$. We use the elementary
matrix-coefficient argument: if $U,D\ne0$ and
$UXA=BXD$ for every compatible full matrix $X$, then
$A=\alpha D$ and $B=\alpha U$ for a scalar $\alpha$.
Indeed, comparison of the coefficient of $X_{ij}$ in entry $(p,q)$
gives $U_{pi}A_{jq}=B_{pi}D_{jq}$. Fix a nonzero entry of $U$,
then one of $D$, to obtain the assertion.

The fixed-coupling cancellation equation is
\[
 WX(\beta CV)=(WN\beta)XV\qquad\text{for every }X.
\]
It follows that $WN\beta=\alpha W$ and $\beta CV=\alpha V$.
For $U=\beta C$, the identity $U-U^2=\beta\nu$ and
$\nu V=0$ imply $\alpha(1-\alpha)V=0$, so $\alpha=0$ or $1$.
At $\alpha=0$, $WN=0$ and $CV=0$; at $\alpha=1$,
$WC=0$ and $NV=0$. These give~\eqref{eq:full-annihilation},
whose converse follows from the unrestricted result already proved.

If $N\nu=0=\nu C$, then $\rk N,\rk C\le m-r$.
Lemma~\ref{lem:rank} gives
$r=\rk N+\rk C-m\le m-2r$, hence $m\ge3r$; the other
orientation is symmetric. For sharpness, take $r$ direct sums of
$N_*=E_{12}+E_{23}$, $C_*=E_{12}+E_{31}$.
Each residue is $E_{12}/2$ with $N_*\nu_*=0=\nu_*C_*$.
The sum pair has $m=3r$, residue of type $[2^r,1^r]$, and
universal cancellation even for directions with off-block entries,
by the criterion just proved. The numerical bound alone is not
sufficient for arbitrary inputs: in $M_3(\C)$, take
$N=E_{12}+E_{33}$ and $C=I-N$. Then $S=I$ and $\nu=E_{12}$,
but $C\nu=\nu C=\nu\ne0$, so neither alternative holds.

Now suppose $\nu^2=0$, $m=4$, $r=2$. Choose a two-column
matrix $Z_0$ with $WZ_0=I_2$. Since $\ker W=\im V$,
the matrix $P=(V,Z_0)$ is invertible, and
\[
 P^{-1}\nu P=\begin{pmatrix}0&I_2\\0&0\end{pmatrix},\qquad
 P^{-1}T'_0P=\begin{pmatrix}H_{11}&H_{12}\\\Gamma&H_{22}\end{pmatrix}.
\]
In the determinant of $P^{-1}T_\varepsilon P$, the two lower
rows are $O(\varepsilon)$. The coefficient of $\varepsilon^2$
uses the two upper-right identity entries and the lower-left
block $\varepsilon\Gamma$, with positive sign. This proves
\eqref{eq:rank-two-determinant}, independently of higher output
coefficients. For $\varepsilon=\zeta^2$ and
$D_\zeta=\operatorname{diag}(I_2,\zeta I_2)$,
\[
 \zeta^{-1}D_\zeta^{-1}P^{-1}T_{\zeta^2}PD_\zeta
 =
 \begin{pmatrix}
  \zeta H_{11}+O(\zeta^3)&I_2+O(\zeta^2)\\
  \Gamma+O(\zeta^2)&\zeta H_{22}+O(\zeta^3)
 \end{pmatrix}.
\]
Each block extends holomorphically to $\zeta=0$.
Its limiting characteristic polynomial is $\det(w^2I_2-\Gamma)$.
Continuity of polynomial roots gives the limiting multiset
$\{\sqrt{\gamma_1},-\sqrt{\gamma_1},
\sqrt{\gamma_2},-\sqrt{\gamma_2}\}$, proving
\eqref{eq:rank-two-leading}. The opening ramification argument
gives the claimed local parametrizations even when the eigenvalues
of $\Gamma$ repeat; no simplicity or stronger remainder is assumed here.

For the nilpotent-input specialization, Lemma~\ref{lem:rank}
forces both input ranks to be three. Thus both adjugates are
nonzero and~\eqref{eq:rank-two-kappa} is a nonempty Zariski-open
condition. The exact determinant factorization gives
$\det T_\varepsilon=\kappa_4(X)\varepsilon^2+O(\varepsilon^3)$,
so comparison with~\eqref{eq:rank-two-determinant} proves the
coefficient identity. The nonzero input determinant derivatives
give punctured invertibility. Alternatively,
$T_\varepsilon^2=O(\varepsilon)$ bounds every root above by a
multiple of $|\varepsilon|^{1/2}$; dividing the determinant
lower bound by the other three upper bounds gives the
corresponding lower bound for each root.
\end{proof}

The square-zero hypothesis in the fixed-coupling assertion is
essential: $N=C=I$ gives $\nu=I/2$ and $L(X,-X)=0$ for every
$X$, but neither ordered alternative holds.
For a restricted space $\mathcal L$ of allowed raw first derivatives,
full cancellation is exactly
\[
 \mathcal L\subseteq\ker\bigl((X,Y)\mapsto\Gamma(X,Y)\bigr).
\]
At an affine fixed-coupling path, all later output coefficients
are fixed by $X$. For general analytic paths with prescribed first
derivative, finer rate questions also involve the variable
higher input coefficients. The scalar identity of
Proposition~\ref{prop:adjugate} and universal full cancellation
are distinct: the direct sums above constrain the full coefficient
at type $[2^r,1^r]$, although for $r\ge2$ this is not type $[m-1,1]$.

\begin{example}[Square-root splitting for a rank-two residue]\label{ex:ranktwo}
In $M_4(\C)$ let
\[
N=E_{12}+E_{23}+E_{34},\qquad C=E_{21}+E_{14}+E_{43},\qquad X=E_{32}+E_{41}.
\]
Both inputs are regular nilpotents of index four. Their sum has determinant $1$, and
\[
(N+C)^{-1}=\begin{pmatrix}
0&1&0&-1\\1&0&-1&0\\0&0&0&1\\0&0&1&0
\end{pmatrix}.
\]
Multiplication gives $\nu=E_{14}+E_{23}$, of type $[2,2]$. For the fixed-coupling path,
\[
T_\varepsilon=\left[\begin{matrix}0 & \varepsilon & 0 & 1\\- \varepsilon & 0 & 1 & 0\\0 & \varepsilon \left(\varepsilon - 1\right) & 0 & \varepsilon\\\varepsilon \left(\varepsilon + 1\right) & 0 & - \varepsilon & 0\end{matrix}\right]
\]
is obtained directly from $(N+\varepsilon X)(N+C)^{-1}(C-\varepsilon X)$; the invariant calculation is
\[
\det(\lambda I-T_\varepsilon)=\lambda^4-\varepsilon^2,
\qquad \det(N+\varepsilon X)=-\varepsilon,\quad \det(C-\varepsilon X)=\varepsilon.
\]
Thus all four output eigenvalues have leading exponent $1/2$, although the limiting residue is square-zero. These identities follow by direct multiplication and determinant expansion. Adjoining a scalar core as in Section~\ref{sec:corners} and choosing the canonical Drazin inverse coordinates gives a limiting nonzero proper core and classical nearby configurations with core $I$. This example shows that the three-dimensional two-thirds bound does not extend to all rank-two square-zero residues. It also distinguishes scalar trace cancellation from vanishing of the full leading coefficient. With
\[
 V=(e_1,e_2),\qquad W=\begin{pmatrix}e_4^T\\ e_3^T\end{pmatrix},
 \qquad \nu=VW,\qquad WV=0,
\]
the derivative $H=T'_0$ satisfies
\[
 WHV=\begin{pmatrix}1&0\\0&-1\end{pmatrix},\qquad
 \tr(\nu H)=\tr(WHV)=0.
\]
The full two-by-two coefficient is nonzero despite its zero trace
for this direction. For an arbitrary affine fixed-coupling direction
$X$ at the same base, direct substitution gives
\[
 \Gamma(X,0)=\operatorname{diag}(X_{41},-X_{32}),\qquad
 \kappa_4(X)=-X_{32}X_{41},\qquad
 \tr(\nu H(X,0))=X_{41}-X_{32}.
\]
Thus trace cancellation here is direction-specific, whereas
$\tr(\nu^2H)=0$ is vacuous.
Proposition~\ref{prop:index} and this example realize, respectively,
the extremal type $[m-1,1]$ in every $m\ge2$ and type $[2,2]$ in
dimension four. Theorem~\ref{thm:nilpotent-realization} gives all
nilpotent types allowed by the rank bound in every dimension.
Proposition~\ref{prop:full-cancellation} shows that the square-root
modulus law here is generic among fixed-coupling directions at
every rank-two square-zero base with nilpotent inputs in dimension four.
\end{example}

\subsection{Exceptional affine directions}

The exceptional bound in Theorem~\ref{thm:generic} can be refined
by a finite calculation when the input path is affine with fixed sum.
The resulting coefficient loci separate this calculation from the
geometric question retained in Problem~\ref{prob:rate}.

\begin{proposition}[A finite coefficient criterion for affine rates]
\label{prop:affine-coefficient-criterion}
Let $N,C\in M_3(\C)$ be nilpotent with invertible sum, put
$\beta=(N+C)^{-1}$, and assume that $\nu=N\beta C$ is nonzero
and nilpotent. For $X\in M_3(\C)$ with $\kappa(X)=0$, set
\[
 \begin{gathered}
 T_\varepsilon=(N+\varepsilon X)\beta(C-\varepsilon X)
 =\nu+\varepsilon H-\varepsilon^2K,\\
 H=X-X\beta N-N\beta X,\qquad K=X\beta X.
 \end{gathered}
\]
Write $\tau=\tr H$, $\sigma=\tr K$, and
\[
 \begin{aligned}
 u_N&=\tr(\adj(N)X),&v_N&=\tr(\adj(X)N),\\
 u_C&=\tr(\adj(C)X),&v_C&=\tr(\adj(X)C),&d&=\det X.
 \end{aligned}
\]
Then $\chi_{T_\varepsilon}(\lambda)
=\lambda^3+A(\varepsilon)\lambda^2+B(\varepsilon)\lambda+D(\varepsilon)$,
where
\begin{equation}\label{eq:affine-characteristic-coefficients}
 \begin{aligned}
 A(\varepsilon)&=-\tau\varepsilon+\sigma\varepsilon^2,\\
 B(\varepsilon)&=
 \left(\frac{\tau^2-\tr(H^2)}2+\tr(\nu K)\right)\varepsilon^2
 +\bigl(\tr(HK)-\tau\sigma\bigr)\varepsilon^3
 +\frac{\sigma^2-\tr(K^2)}2\varepsilon^4,\\
 D(\varepsilon)&=\det(\beta)\varepsilon^2
 (u_N+\varepsilon v_N+\varepsilon^2d)
 (u_C-\varepsilon v_C+\varepsilon^2d).
 \end{aligned}
\end{equation}
Give an identically zero polynomial order $\infty$. Then
\[
 \operatorname{ord}_0 A\in\{1,2,\infty\},\qquad
 \operatorname{ord}_0 B\in\{2,3,4,\infty\},\qquad
 \operatorname{ord}_0 D\in\{3,4,5,6,\infty\}.
\]
The negative slopes of the lower convex hull of the finite points
\[
 (3,0),\quad (2,\operatorname{ord}_0 A),\quad
 (1,\operatorname{ord}_0 B),\quad (0,\operatorname{ord}_0 D)
\]
are the orders of the branches that are not identically zero, with
multiplicities equal to the horizontal edge lengths. If $k$ is the
smallest power of $\lambda$ whose coefficient is not identically
zero in $\varepsilon$, exactly $k$ branches are identically zero.

For a specified complex linear space $L$ with $\kappa|_L=0$, the
triples of coefficient orders define a finite partition of $L$
into locally closed loci. There are at most $60$ formal valuation triples before
dependencies are imposed. On a nonempty Zariski-open dense subset
of $L$, each coefficient order is the least power of $\varepsilon$ whose
coefficient is not identically zero on $L$, or $\infty$ if none exists.
\end{proposition}
\begin{proof}
Put $a_\varepsilon=N+\varepsilon X$. Expanding
$a_\varepsilon-a_\varepsilon\beta a_\varepsilon$ in the displayed
order gives the formula for $T_\varepsilon$. By Theorem~\ref{thm:mechanism},
$\nu^2=0$ and $\tr(\nu H)=0$. Substitution in
\[
 A=-\tr T_\varepsilon,\qquad
 B=\tfrac12\bigl((\tr T_\varepsilon)^2-\tr(T_\varepsilon^2)\bigr)
\]
gives the displayed expressions. Determinant multilinearity gives
\[
 \begin{aligned}
 \det(N+\varepsilon X)
   &=\varepsilon(u_N+\varepsilon v_N+\varepsilon^2d),\\
 \det(C-\varepsilon X)
   &=-\varepsilon(u_C-\varepsilon v_C+\varepsilon^2d).
 \end{aligned}
\]
Since $D=-\det T_\varepsilon$, the determinant factorization
\eqref{eq:determinant-parallel-path} proves the formula for $D$.
Its quadratic coefficient is $-\kappa(X)=0$, giving the stated
order ranges.

Use the ramified factorization established at the start of this
section. After $\varepsilon=\zeta^r$, discard the identically zero
branches and write each remaining branch as
$c_j\zeta^{n_j}+\cdots$, with $c_j\ne0$ and the $n_j$ in
nondecreasing order. Each characteristic coefficient is an elementary
symmetric sum of branch products. Its order is at least the sum
of the corresponding smallest $n_j$. At each boundary between groups
of equal orders, the least-order product is unique and nonzero.
These boundary points therefore lie on the lower convex hull;
cancellation within a group can only raise intermediate points.
Dividing the orders by $r$ proves the slope rule, including repeated
leading roots. Factoring out the largest common power $\lambda^k$
proves the assertion about permanent zero branches.

All displayed coefficients are polynomial functions of $X$ for
the fixed base. A coefficient order is specified by vanishing of its preceding
coefficients and nonvanishing of its first coefficient, or vanishing
of them all for order $\infty$. This gives locally closed loci and
at most $3\cdot4\cdot5=60$ formal triples. The first coefficient
of each polynomial that is not identically zero on $L$ can be
simultaneously nonzero: their product is nonzero in the domain
$\C[L]$. Its nonvanishing set proves the generic assertion.
\end{proof}

The $60$ formal triples need not all occur, and distinct triples
can have the same lower hull. The generic assertion is relative to
$L$; for a reducible algebraic direction set contained in
$\{\kappa=0\}$ it applies separately to its irreducible components.

\begin{corollary}[Sharp bounds on exceptional affine gaps]
\label{cor:affine-sharp-gap}
Under the hypotheses and notation of
Proposition~\ref{prop:affine-coefficient-criterion}, let $k$ be the
number of identically zero eigenvalue branches, counted with
multiplicity. Every other branch has order at least $1$ and at most
$5/2$ if $k=0$, at most $3$ if $k=1$, and at most $2$ if $k=2$.
Thus, if $k\le2$, there are constants $c,\delta>0$, depending on
the fixed path, such that
\begin{equation}\label{eq:affine-sharp-gap}
 \gap(T_\varepsilon)\ge c|\varepsilon|^3,
 \qquad 0<|\varepsilon|<\delta.
\end{equation}
The exponent $3$ is optimal over this class, and $5/2$ is optimal
for the largest branch order in the subclass $k=0$.
If $k=3$, the gap is $+\infty$.
Moreover, when $D\not\equiv0$, the coefficients satisfy
\begin{equation}\label{eq:affine-order-constraint}
 \operatorname{ord}_0 B\ge\operatorname{ord}_0 D-2,
\end{equation}
with the convention $\operatorname{ord}_0 0=\infty$.
\end{corollary}
\begin{proof}
Theorem~\ref{thm:generic} gives order at least $1$ for every moving
branch. If $k=1$, the two moving orders sum to
$\operatorname{ord}_0 B\le4$, so each is at most $3$. If $k=2$,
the sole moving order is $\operatorname{ord}_0 A\le2$.
If $k=3$, there is no nonzero eigenvalue.

Suppose now that $k=0$, equivalently $D\not\equiv0$.
Put $U_\varepsilon=N+\varepsilon X$ and
$V_\varepsilon=C-\varepsilon X$. The determinant factorization
shows that both input determinants are nonzero polynomials.
Each vanishes at zero and has degree at most three. On a
sufficiently small punctured disk both inputs and $T_\varepsilon$
are therefore invertible, and ordered inversion gives
\[
 T_\varepsilon^{-1}
 =V_\varepsilon^{-1}(U_\varepsilon+V_\varepsilon)U_\varepsilon^{-1}
 =V_\varepsilon^{-1}+U_\varepsilon^{-1}.
\]
For a nilpotent $3\times3$ matrix $M$, $\tr(\adj M)=0$.
Consequently $\tr(\adj U_\varepsilon)$ and
$\tr(\adj V_\varepsilon)$ are $O(|\varepsilon|)$.
Dividing by their respective determinants, each of order at most
three, yields
\[
 \tr(T_\varepsilon^{-1})=O(|\varepsilon|^{-2}).
\]
After a common ramification, order the three branch orders as
$1\le s_1\le s_2\le s_3$. They sum to
$\operatorname{ord}_0 D\le6$. If $s_3>s_2$, the reciprocal of
the third branch is the unique term of largest pole order in
$\tr(T_\varepsilon^{-1})$, so its nonzero leading term cannot
cancel; hence $s_3\le2$. If $s_3=s_2$, then
$2s_3+s_1\le6$, whence $s_3\le5/2$.
The identity
$B=\det(T_\varepsilon)\tr(T_\varepsilon^{-1})$
also proves~\eqref{eq:affine-order-constraint}.

Every moving branch has a convergent ramified expansion with a
nonzero leading coefficient. Taking a common sufficiently small
disk and the minimum of finitely many positive lower constants
therefore proves~\eqref{eq:affine-sharp-gap}, since all moving
orders are at most $3$. This holds for complex parameters.
The canonical example in Remark~\ref{rem:census} has orders
$[1,5/2,5/2]$ and no permanent zero, proving sharpness for $k=0$.
Example~\ref{ex:affine-cubic-gap} below has gap asymptotic to
$4|\varepsilon|^3$, proving optimality of the cubic bound over
the whole exceptional affine class.
\end{proof}

These conclusions concern affine fixed-sum input paths. General
analytic input paths with unrestricted higher coefficients are
outside their scope.

\begin{example}[An exceptional affine path with cubic gap]
\label{ex:affine-cubic-gap}
Take
\[
 N=E_{12}+E_{23},\qquad C=E_{12}+E_{31},\qquad
 X=\operatorname{diag}(1,2,0)+E_{13}.
\]
Both base inputs have index three, $\det(N+C)=2$ and
$\nu=E_{12}/2$. Here $\kappa(X)=0$, and direct multiplication gives
\[
 \begin{gathered}
 T_\varepsilon=
 \begin{pmatrix}
 \varepsilon/2&1/2-2\varepsilon^2&-\varepsilon/2\\
 -\varepsilon^2&-\varepsilon&-\varepsilon^2\\0&0&0
 \end{pmatrix},\\
 \chi_{T_\varepsilon}(\lambda)
 =\lambda(\lambda^2+\varepsilon\lambda/2-2\varepsilon^4).
 \end{gathered}
\]
There is one permanent zero branch. The other two are
\[
 \frac{\varepsilon}{4}\bigl(-1\pm\sqrt{1+32\varepsilon^2}\bigr),
\]
using the analytic square root equal to $1$ at zero. They have expansions
$-\varepsilon/2-4\varepsilon^3+O(\varepsilon^5)$ and
$4\varepsilon^3+O(\varepsilon^5)$. Consequently
\[
 \max_j|\lambda_j(\varepsilon)|\sim|\varepsilon|/2,
 \qquad \gap(T_\varepsilon)\sim4|\varepsilon|^3.
\]
This attains the cubic bound in Corollary~\ref{cor:affine-sharp-gap}
over the whole exceptional affine class, as well as the upper
order $3$ in the subclass with one permanent zero.
The input determinants are $\det(N+\varepsilon X)=0$ and
$\det(C-\varepsilon X)=-2\varepsilon^2$. The path has fixed
invertible sum but does not retain nilpotence of its inputs.
The output remains singular; its gap excludes the permanent zero.
\end{example}

\section{Coefficient algebras and Hua configurations}\label{sec:corners}

\subsection{Coefficient-ring realization and Banach spectra}
We now allow coefficients in a unital ring. The construction gives
inputs of index exactly three and, in the Banach setting, prescribed
spectra after matrix amplification.

\begin{theorem}[Realization over a coefficient algebra]\label{thm:realization}
Let $\K$ be a nonzero unital ring in which $2=1+1$ is invertible,
and let $Q\in\K$. In $M_3(\K)$ put
\begin{equation}\label{eq:Q-inputs}
 N=\begin{pmatrix}0&1&0\\0&0&1\\0&0&0\end{pmatrix},\qquad
 C_Q=\begin{pmatrix}-Q&1&0\\-Q^2&Q&0\\1&0&0\end{pmatrix}.
\end{equation}
Both inputs are nilpotent of index exactly three. Their sum is invertible, with
\begin{equation}\label{eq:Q-inverse}
 (N+C_Q)^{-1}=
 \begin{pmatrix}0&0&1\\\tfrac12&0&\tfrac Q2\\-\tfrac Q2&1&\tfrac{Q^2}2\end{pmatrix},
\end{equation}
and their residue is
\begin{equation}\label{eq:Q-residue}
 R_Q=\Par(N,C_Q)=\frac12
 \begin{pmatrix}0&1&0\\0&Q&0\\0&0&0\end{pmatrix}.
\end{equation}
Every $R_Q$ has the same square-zero reflexive inverse $G=2E_{21}$:
\[
 R_QGR_Q=R_Q,\qquad GR_QG=G,\qquad G^2=0.
\]
This inverse does not commute with $R_Q$.
If $\K$ is a complex unital Banach algebra, the spectrum of $R_Q$
in $M_3(\K)$ is
\begin{equation}\label{eq:Q-spectrum}
 \sigma(R_Q)=\{0\}\cup\tfrac12\sigma_{\K}(Q).
\end{equation}
\end{theorem}
\begin{proof}
Here $1$ in a matrix entry denotes the unit of $\K$. We have $N^2=E_{13}\ne0$, $N^3=0$, and
\[
 C_Q^2=\begin{pmatrix}0&0&0\\0&0&0\\-Q&1&0\end{pmatrix}\ne0,
 \qquad C_Q^3=0.
\]
The entries $1$ are nonzero because $\K$ is nonzero. Multiplying the matrix in~\eqref{eq:Q-inverse} by
\[
 N+C_Q=\begin{pmatrix}-Q&2&0\\-Q^2&Q&1\\1&0&0\end{pmatrix}
\]
in either order gives $I_3$. Multiplication by $N$ on the left and $C_Q$ on the right gives~\eqref{eq:Q-residue}. These identities use only powers of $Q$ and scalar entries, so no commutativity assumption on $\K$ is needed.

For $G=2E_{21}$, direct multiplication gives
\[
 R_QG=E_{11}+QE_{21},\qquad GR_Q=E_{22}.
\]
These ordered products give the stated reflexive identities, and
$G^2=0$. The products differ because $\K$ is nonzero, proving
noncommutation. No commutativity of $\K$ is required.

Now suppose that $\K$ is a complex unital Banach algebra. For $\lambda\ne0$ there is the exact factorization
\begin{equation}\label{eq:Q-factorization}
 \lambda I_3-R_Q=
 \operatorname{diag}(\lambda,\lambda-Q/2,\lambda)
 \left(I_3-\frac1{2\lambda}E_{12}\right).
\end{equation}
The second factor is a unit, with inverse $I_3+E_{12}/(2\lambda)$. A diagonal matrix is invertible exactly when its diagonal entries are invertible: the corresponding diagonal entries of any two-sided inverse provide both inverses of each entry. Thus, for $\lambda\ne0$, the left side of~\eqref{eq:Q-factorization} is invertible exactly when $\lambda-Q/2$ is invertible in $\K$.

At zero, $R_QE_{11}=0$ and $E_{11}\ne0$, so $R_Q$ cannot be invertible. This proves~\eqref{eq:Q-spectrum}, including its added zero. The argument uses the displayed factorization, not a converse inferred from upper triangularity with operator entries.
\end{proof}

The construction realizes a scaled spectrum after matrix amplification. It does not prescribe an arbitrary output matrix in the same fixed dimension as $Q$.

In dimension three, Theorem~\ref{thm:complete-realization} says that
all residues are precisely the matrices of rank at most one; their
nilpotent members are zero and the nonzero rank-one nilpotents.
The scalar specialization of Theorem~\ref{thm:realization} gives a
short explicit realization of every nonzero rank-one similarity class:
\[
 R_x=\tfrac12E_{12}+\tfrac x2E_{22},\qquad
 R_x^2=\tfrac x2R_x\qquad(x\in\C).
\]
Here $C_x=(I+xE_{21})C_0(I-xE_{21})$ and
$\sigma(R_x)=\{0,x/2\}$. Both input indices are three for every $x$.
At $x=0$, $R_x$ is
nonzero nilpotent; at $x=2\lambda$ it has prescribed nonzero trace
$\lambda$. Indeed every rank-one matrix $R=uv^T$ satisfies
$R^2=(\tr R)R$. If $\tr R=\lambda\ne0$, its minimal polynomial
divides $z(z-\lambda)$, so it is similar to
$\operatorname{diag}(\lambda,0,0)$. If its trace is zero, its
Jordan type is $[2,1]$, so it is similar to $E_{12}$.
Simultaneous conjugation of the inputs conjugates their residue.
Zero is also obtained directly from $N=E_{12}$ and
$C=E_{23}+E_{31}$: their sum is an invertible cycle and their
parallel sum is zero.

In a complex Banach algebra, $x$ is quasinilpotent when $\rho(x)=0$. Its generalized Drazin inverse is the unique $y$ satisfying $xy=yx$, $yxy=y$, and quasinilpotence of $x-x^2y$. It exists exactly when zero is not an accumulation point of $\sigma(x)$~\cite[Theorems 3.1, 4.2 and 4.3]{Koliha96}. A finite Drazin inverse has a nilpotent remainder instead. A quasinilpotent element has generalized Drazin inverse zero, and has a finite Drazin inverse exactly when it is nilpotent. These are Banach-algebra assertions; the spectrum in~\eqref{eq:Q-spectrum} is taken in the specified unital algebra.

The following boundary instances also appear in~\cite[Example~7.4]{HuaI};
we retain their verification from the general construction.
\begin{example}[Quasinilpotence and accumulation at zero]\label{ex:Banach}
Two choices of $Q$ in Theorem~\ref{thm:realization} separate the possible output regimes.

First take the Volterra operator $V$ on $C[0,1]$, defined by $(Vf)(r)=\int_0^r f(u)\,du$, in the Banach algebra of bounded operators on that space. Repeated integration gives
\[
 (V^kf)(r)=\frac1{(k-1)!}\int_0^r(r-u)^{k-1}f(u)\,du
 \quad(k\ge1).
\]
Hence $\|V^k\|=1/k!$, with equality on the constant function one, and $V^k\ne0$ for every $k$. Direct induction gives
\[
 R_V^j=2^{-j}\begin{pmatrix}0&V^{j-1}&0\\0&V^j&0\\0&0&0\end{pmatrix}
 \quad(j\ge1),\qquad V^0=I.
\]
Every such power is nonzero, whereas its root norm tends to zero by the factorial estimate. Thus $R_V$ is quasinilpotent and not nilpotent. Its generalized Drazin inverse is zero, but it has no finite Drazin inverse, although both original inputs have index three.

Second let $Q$ act on $\ell^2$ diagonally by $Qe_n=(2/n)e_n$. Its spectrum is $\{0\}\cup\{2/n:n\ge1\}$. Indeed each $2/n$ is an eigenvalue; $Q$ is not invertible because $\|Qe_n\|\to0$; and outside this closed set the diagonal reciprocals of $\lambda-Q$ are uniformly bounded and define its inverse. Equation~\eqref{eq:Q-spectrum} now gives
\[
 \sigma(R_Q)=\{0\}\cup\{1/n:n\ge1\}.
\]
Zero is an accumulation point, so $R_Q$ has no generalized Drazin inverse. Quasinilpotent nonnilpotent residues have no finite-dimensional matrix instance, since a complex matrix with spectrum $\{0\}$ is nilpotent by Cayley--Hamilton.
\end{example}

Both Banach examples above have the bounded square-zero reflexive
inverse $G=2E_{21}$ from Theorem~\ref{thm:realization}.
The Volterra residue has no finite Drazin inverse, while the diagonal
example has no Banach generalized Drazin inverse.

\subsection{Common-core completion and finite matrix outputs}
For the connection with Hua inversion, selected outer-inverse data satisfy
\[
 a'aa'=a',\quad aa'=a'a=p,\qquad
 c'cc'=c',\quad cc'=c'c=r.
\]
Each selected outer inverse commutes with its own input, and $p,r$
are the resulting core idempotents. The element intertwining
condition is $pb=br$. An augmented tuple includes these selected
inverse coordinates as well as the raw inputs.

For complex matrices, also require $b\in GL_n(\C)$, $a+c=b^{-1}$
and nilpotent input remainders $n_a=a(I-p)$, $n_c=c(I-r)$.
Together with $pb=br$, these conditions define the finite-Drazin
Hua configurations used below; then $a'=a^{\Dz}$ and $c'=c^{\Dz}$.
This is the convention of~\cite[Definition~2.1 and Notation~2.3]{HuaI}.
The following common-core construction is also valid in a unital ring.

The preceding constructions can be completed by an invertible core. Let $\A$ be a unital ring, let $e\in\A$ be an idempotent, and put $f=1-e$. Take units $A_0,C_0,A_0+C_0$ in $e\A e$ and nilpotent elements $N,C\in f\A f$ whose sum is a unit of that corner. All inverses in the next display are taken in their indicated corners. Set
\begin{align}\label{eq:completion}
 a&=A_0+N,&c&=C_0+C,&
 b&=(A_0+C_0)^{-1}_{e\A e}+\beta,\\
 a'&=(A_0)^{-1}_{e\A e},&c'&=(C_0)^{-1}_{e\A e},&
 \beta&=(N+C)^{-1}_{f\A f}.\notag
\end{align}
Products between the two corners vanish. Thus $a+c=b^{-1}$, $a'aa'=a'$, $c'cc'=c'$, and $aa'=a'a=cc'=c'c=e$. Thus $a'$ commutes with $a$ and $c'$ with $c$, each is an outer inverse of its input, and $p=r=e$. Since $eb=be$, the element intertwining condition $pb=br$ holds. The input remainders are
\[
 n_a=a-aa'a=N,\qquad n_c=c-cc'c=C.
\]
The Hua output and its prescribed core inverse are
\begin{equation}\label{eq:core-output}
 t=a-aba=T_{\rm core}+\nu,\qquad
 T_{\rm core}=A_0(A_0+C_0)^{-1}_{e\A e}C_0,\qquad s=a'+c'.
\end{equation}
Indeed $T_{\rm core}^{-1}=C_0^{-1}(A_0+C_0)A_0^{-1}=C_0^{-1}+A_0^{-1}$, and the complementary block is~\eqref{eq:symmetry}. In particular $ts=st=e$, $sts=s$, and $t-t^2s=\nu$. This verifies the completion directly. We allow $e=0$; the zero corner has its own unit zero. In the Banach setting its corner spectrum is empty, with distance $+\infty$ from zero.

For a scalar core choose $\alpha,\gamma\in\C\setminus\{0\}$ with $\alpha+\gamma\ne0$. In $\C\oplus B$, any nilpotent pair $N,C\in B$ with invertible sum gives
\begin{align}\label{eq:scalar-completion}
 a&=\alpha\oplus N,&c&=\gamma\oplus C,&
 b&=(\alpha+\gamma)^{-1}\oplus(N+C)^{-1},\\
 e&=1\oplus0,&
 t&=\frac{\alpha\gamma}{\alpha+\gamma}\oplus\nu,&
 s&=(\alpha^{-1}+\gamma^{-1})\oplus0.\notag
\end{align}
For coupled data with common core idempotent $e$, we use
\emph{classical} for $e=1$ and \emph{nondegenerate} for
$0\ne e\ne1$. The scalar-core completion above is nondegenerate
when $B\ne0$. For nilpotent $\nu$, $s=t^{\Dz}$; when $B\ne0$, the index of $t$ is the nilpotency index of $\nu$. If $B=0$, then $t$ is a unit and its Drazin index is zero. In the Banach setting, $s=t^{\gd}$ exactly when $\nu$ is quasinilpotent; the mere existence of an output generalized Drazin inverse is a different question. These conclusions follow directly from $t-t^2s=\nu$ and the displayed block decomposition, or from~\cite[Theorems~5.2 and~5.4]{HuaI}.

The unit scalar-core instances $\alpha=\gamma=1$ below also
appear in~\cite[Examples~7.1--7.3]{HuaI}; we retain the full parameter families.
\begin{example}[Scalar-core families]\label{ex:scalar-families}
Taking $B=M_2(\C)$, $N=E_{12}$ and $C=E_{21}$ in~\eqref{eq:scalar-completion} gives
\[
 a=\alpha E_{11}+E_{23},\quad c=\gamma E_{11}+E_{32},\quad
 t=\frac{\alpha\gamma}{\alpha+\gamma}E_{11}.
\]
The completed inputs $a,c$ both have finite Drazin index two,
$\nu=0$, and $t^{\#}=s$. The remainders do not commute; also $ab e_2=e_2$ while $ba e_2=0$, so the inputs do not commute with the coupling unit.

Taking instead the $Q=0$ pair from~\eqref{eq:Q-inputs} in $B=M_3(\C)$ gives, in $M_4(\C)$,
\[
 a=\alpha E_{11}+E_{23}+E_{34},\qquad
 c=\gamma E_{11}+E_{23}+E_{42},
\]
\[
 b=(\alpha+\gamma)^{-1}E_{11}+\tfrac12E_{32}+E_{43}+E_{24},
 \qquad
 t=\frac{\alpha\gamma}{\alpha+\gamma}E_{11}+\tfrac12E_{23}.
\]
Here $a,c$ both have finite Drazin index three, while $t$ has
finite Drazin index two and $t^{\Dz}=s$; in particular, $t$ has
no group inverse. These formulas also hold over $\mathbb Q$ for rational $\alpha,\gamma$ satisfying the same nonvanishing conditions. At $\alpha=\gamma=1$, the resolvent is
\[
 (\lambda-t)^{-1}=(\lambda-\tfrac12)^{-1}E_{11}
 +\lambda^{-1}(I-E_{11})+\tfrac12\lambda^{-2}E_{23},
\]
with a double pole at zero. Omitting the scalar core gives the fully degenerate ($e=0$) $M_3(\C)$ pair with $\nu=E_{12}/2$. Finally, adjoining the scalar core to the two pairs in Example~\ref{ex:Banach} preserves their residual spectral behavior in nondegenerate configurations.
\end{example}

\begin{remark}[Dimension and noncommutativity]\label{rem:dimensions}
Two commuting nilpotents cannot have invertible sum in a nonzero corner, because their sum is nilpotent by the binomial expansion. More generally the common-idempotent ring configurations of~\cite[Proposition 6.4]{HuaI} have no nonzero commutative complementary corner when the input remainders are Harte-quasinilpotent. Here that term means that $1+qy$ is invertible for every $y$ commuting with $q$. To see the obstruction directly in the corner, if $u=q_1+q_2$ is a unit and the summands commute, then $u^{-1}$ commutes with them. Quasinilpotence makes $1-q_1u^{-1}=q_2u^{-1}$ a unit, and hence makes $q_2$ a unit, impossible in a nonzero ring. Compression preserves this quasinilpotence for a remainder $q$
commuting with $f$, as in~\cite[Lemma~6.1]{HuaI}. Indeed, if
$y\in f\A f$ commutes with $fqf$, it commutes with $q$.
The unit $1+qy$ commutes with $f$, so its $f$-corner
$f+(fqf)y$ is a unit in $f\A f$.

Consequently, for complex matrix configurations with quasinilpotent input remainders, the common core cannot have corank one. In $M_2(\C)$ it is either $e=0$ or $e=I$; if $e=0$, the outer equations force $a'=c'=0$. Coexistence of core and complement requires total dimension at least three, and the first family in Example~\ref{ex:scalar-families} attains that dimension.

Every quasinilpotent matrix in $M_2(\C)$ is square-zero. Proposition~\ref{prop:square-zero} therefore forces a zero residue in a two-dimensional complementary corner. A nonzero nilpotent residue needs complementary dimension at least three, and a nondegenerate such configuration needs total dimension at least four, attained by the second family above. The assumption $0\ne e\ne I$ is essential for the total-dimension bound: the fully degenerate pair already has nonzero residue in dimension three.
\end{remark}

\begin{corollary}[Finite-dimensional Hua outputs]
\label{cor:hua-output-classification}
Let $n\ge1$ and $T\in M_n(\C)$.
The matrix $T$ is the output of a finite-Drazin Hua configuration
if and only if $\dim\ker T\ne1$.
Here the common core may be $0$ or $I$; the assertion ranges over
all such configurations, without prescribing a core.
\end{corollary}
\begin{proof}
The common-core reduction is proved in
\cite[Theorems~2.5 and~3.2]{HuaI}; we give a direct matrix argument.
Set $d=bcb^{-1}$ and $d'=bc'b^{-1}$.
Since $pb=br$, one has $dd'=d'd=brb^{-1}=p$.
Relative to $p\C^n\oplus(I-p)\C^n$, write
$a=a_+\oplus a_-$ and $d=d_+\oplus d_-$.
The plus blocks are invertible, with inverses supplied by $a'$ and
$d'$, and the minus blocks are nilpotent. For $d_-$ this follows
from $d(I-p)=bc(I-r)b^{-1}$, without assuming $bp=pb$.
The off-diagonal blocks of $ba+db=b(a+c)=I$ satisfy
\[
 d_+b_{+-}+b_{+-}a_-=0,\qquad
 b_{-+}a_++d_-b_{-+}=0.
\]
Ordered iteration gives
\[
 b_{+-}=(-d_+^{-1})^j b_{+-}a_-^j=0,\qquad
 b_{-+}=(-d_-)^k b_{-+}a_+^{-k}=0
\]
for sufficiently large $j,k$. Hence $bp=pb$, and $pb=br$ now
implies $p=r=e$. Thus $a,c,b$ are block diagonal for the common
core. The two block sums are invertible because $a+c=b^{-1}$.
Writing $m=\dim\ker e$, $N=a_-$ and $C=c_-$ gives, after the
corresponding simultaneous similarity,
\[
 \begin{gathered}
 t=T_{\rm core}\oplus\nu,\qquad \nu=N(N+C)^{-1}C,\\
 T_{\rm core}=a_+(a_++c_+)^{-1}c_+\in GL_{n-m}(\C).
 \end{gathered}
\]
Indeed $T_{\rm core}^{-1}=c_+^{-1}+a_+^{-1}$, and $N,C$ are
nilpotent. The core or complementary summand is allowed to be absent.
If $m=0$, then $t$ is invertible.
If $m>0$, Lemma~\ref{lem:rank} gives
\[
 \dim\ker t=\dim\ker\nu
          =\dim\ker N+\dim\ker C\ge2.
\]
In particular, a one-dimensional complementary summand cannot occur.

Conversely, if $T$ is invertible, take
\[
 a=c=2T,\qquad b=(4T)^{-1},\qquad
 a'=c'=(2T)^{-1}.
\]
These data have $p=r=I$, zero input remainders and $a-aba=T$.
If $\dim\ker T\ge2$, Theorem~\ref{thm:complete-realization}
gives nilpotent $N,C$ with invertible sum and
$T=N(N+C)^{-1}C$. Take
$a=N$, $c=C$, $b=(N+C)^{-1}$ and $a'=c'=0$.
These data have $p=r=0$, satisfy all the displayed conditions,
and have output $T$.
\end{proof}

\begin{corollary}[Generic failure of the zero-core candidate]
\label{cor:generic-candidate-failure}
Let $m\ge3$, and put
\[
 \mathcal U_m=\{(N,C)\in M_m(\C)^2:N^m=C^m=0,
                                   \ \det(N+C)\ne0\}.
\]
The pairs for which $\Par(N,C)$ is not nilpotent form a nonempty
Zariski-open dense subset of $\mathcal U_m$.
For the zero-core Hua configurations $a=N$, $c=C$,
$b=(N+C)^{-1}$ and $a'=c'=0$, the candidate $s=0$ therefore
fails to equal the output's finite Drazin inverse on that subset.
The output itself always has a finite Drazin inverse.
\end{corollary}
\begin{proof}
Let $\mathcal N_m=\{A:A^m=0\}$, and let $\mathcal O_m$ be
the similarity orbit of the regular nilpotent block $J_m$.
In any nilpotent Jordan form, replace the missing superdiagonal
entries between successive blocks by a parameter $z$.
For $z\ne0$ all chain weights are nonzero, so the resulting
matrix is similar to $J_m$; at $z=0$ it is the original form.
Conjugating and using the same parameter for both members of a
pair gives
\[
 \mathcal N_m\times\mathcal N_m
   =\overline{\mathcal O_m\times\mathcal O_m}.
\]
Indeed, a polynomial vanishing on the product orbit vanishes on
each such path for all nonzero $z$, and hence at $z=0$.
The product orbit is the image of the irreducible principal open
set $GL_m(\C)\times GL_m(\C)$ in affine space under
$(g,h)\mapsto(gJ_mg^{-1},hJ_mh^{-1})$.
Its closure is therefore irreducible.
Theorem~\ref{thm:complete-realization} makes $\mathcal U_m$
a nonempty open subset, so it too is irreducible.

On $\mathcal U_m$, nilpotence of the parallel sum is equivalent to
\[
 \bigl(N\adj(N+C)C\bigr)^m=0,
\]
which is a polynomial condition after the nonzero determinant of
$N+C$ has been cleared. This relatively closed subset is proper:
the rank-one nonnilpotent target $E_{11}$ is realizable when
$m\ge3$. Its complement is consequently nonempty, open and dense.
Finally, \cite[Theorem~5.2(c)]{HuaI} identifies success of the
prescribed candidate exactly with nilpotence of the residue.
Every finite complex matrix has a finite Drazin inverse by its
Fitting decomposition, independently of that candidate condition.
\end{proof}

The density in Corollary~\ref{cor:generic-candidate-failure} is in
$\mathcal U_m$ as a whole; no density in every fixed input Jordan
stratum, fixed-sum fibre or prescribed-core class is asserted.

\subsection{Affine corner laws and separated clusters}

For a complex unital Banach algebra, use the same definition of
$\gap(x)$ with its algebra spectrum, and write $\rho(x)$ for the
spectral radius. Spectra and spectral radii of complementary elements
are taken in the indicated corner algebra.

The algebraic identities in Propositions~\ref{prop:deglaw}
and~\ref{prop:compression} hold for common-core data in a unital
ring; Proposition~\ref{prop:quadratic} additionally assumes a
complex unital Banach algebra. We require the units $A_0,C_0,A_0+C_0$ in $e\A e$ and $N+C$ in $f\A f$, but do not require $N,C$ to be nilpotent. Put $n_a=N$, $n_c=C$, $\beta=(N+C)^{-1}_{f\A f}$, and let $T_{\rm core}$ be the core in~\eqref{eq:core-output}. The corner equations used below follow from these units alone. In the ring setting, a derivative of these polynomial expressions means the coefficient of a central formal parameter; in matrix or Banach settings it is the usual derivative.

\begin{proposition}[Canonical scalar-shift law]\label{prop:deglaw}
For these common-core data, put $f=1-e$. For a central scalar $\varepsilon$, let
$a_\varepsilon=a+\varepsilon f$ and $c_\varepsilon=c-\varepsilon f$, with $b$ fixed. Then $t_\varepsilon=a_\varepsilon-a_\varepsilon ba_\varepsilon$ has zero off-corner blocks, unchanged core block $T_{\rm core}$, and complementary block
\[
f t_\varepsilon f=\nu+\varepsilon(n_c\beta-\beta n_a)-\varepsilon^2\beta.
\]
Moreover $n_c\beta=\beta n_a$ if and only if $n_a^2=n_c^2$.
\end{proposition}
\begin{proof}
All limiting data commute with $e$, hence so do the perturbed inputs. In the complementary corner, whose unit is $f$,
\[
(n_a+\varepsilon f)\beta(n_c-\varepsilon f)
=\nu+\varepsilon(\beta n_c-n_a\beta)-\varepsilon^2\beta.
\]
The two forms of the linear coefficient agree because $\beta(n_a+n_c)=(n_a+n_c)\beta=f$. Multiplying $n_c\beta=\beta n_a$ on both sides by $n_a+n_c$ gives $(n_a+n_c)n_c=n_a(n_a+n_c)$, which is equivalent to $n_c^2=n_a^2$.
\end{proof}

\begin{proposition}[A quadratic gap family with zero limiting residue]\label{prop:quadratic}
In the setting of Proposition~\ref{prop:deglaw}, let the ambient algebra be a complex unital Banach algebra $\A$, with common idempotent $e$ and nonzero complementary corner $B=f\A f$, $f=1-e$. Suppose $\nu=0$ and $n_a^2=n_c^2$. Then $n_a^2=n_c^2=0$. For $\varepsilon\in\C\setminus\{0\}$, both perturbed inputs are invertible and
\[
f t_\varepsilon f=-\varepsilon^2\beta,
\qquad f s_\varepsilon f=-\varepsilon^{-2}(n_a+n_c),
\]
where $s_\varepsilon=a_\varepsilon^{-1}+c_\varepsilon^{-1}=t_\varepsilon^{-1}$. For the full element,
\[
\gap(t_\varepsilon)
=\min\left\{d(0,\sigma_{e\A e}(T_{\rm core})),
\frac{|\varepsilon|^2}{\rho_B(n_a+n_c)}\right\}.
\]
The second term is the gap for sufficiently small nonzero $\varepsilon$. The zero-core convention makes the first term $+\infty$ if $e=0$.
\end{proposition}
\begin{proof}
Write $N=n_a,C=n_c$ in the complementary corner. Since $\nu=0$, $N=N\beta N$. The equality of the squares gives $\beta N=C\beta$, so
\[
N=N\beta N=NC\beta,\qquad N(N+C)=NC.
\]
Thus $N^2=0$, and $C^2=0$ follows. The shifted corner inputs have inverses
$\varepsilon^{-1}f-\varepsilon^{-2}N$ and
$-\varepsilon^{-1}f-\varepsilon^{-2}C$. Their core blocks are unchanged and invertible. Proposition \ref{prop:deglaw} gives the output block, and summing the input inverses gives the inverse formula. Spectral splitting and $\sigma_B(\beta)=\{\lambda^{-1}:\lambda\in\sigma_B(N+C)\}$ prove the minimum formula. Since $N+C$ is invertible in a nonzero Banach algebra, its spectral radius is positive.
\end{proof}

\begin{proposition}[First-order compression and its quantifier]\label{prop:compression}
In a fixed complementary corner with unit $f$, put $P_\ell=n_a\beta$, $P_r=\beta n_a$. For direction $X$ and fixed $b$ the derivative of the complementary output is
\[
DF_{(n_a,\beta)}[X,0]=X(f-P_r)-P_\ell X.
\]
For one direction it equals $[X,P_\ell]$ if and only if $X(f-P_r-P_\ell)=0$. It equals $[X,P_\ell]$ for every direction $X$ if and only if $P_\ell+P_r=f$.
\end{proposition}
\begin{proof}
Differentiate $(n_a+\tau X)\beta(n_c-\tau X)$ and subtract $[X,P_\ell]$. The difference is $X(f-P_r-P_\ell)$. For the universal implication choose $X=f$.
\end{proof}

The universal commutator identity is equivalent to
$C\beta=\beta N$, hence to $N^2=C^2$ by
Proposition~\ref{prop:deglaw}. In the nonzero nilpotent
three-dimensional normal form, $N^2=E_{13}$ and $C^2=uvE_{32}$
are distinct. Scalar cancellation there therefore comes from the
ordered annihilation mechanism, while a commutator formula for one
chosen direction retains the separate quantifier of
Proposition~\ref{prop:compression}.

\begin{corollary}[A separated three-dimensional cluster]\label{cor:cluster}
Let $a_0,b_0,c_0$ be limiting matrix data from a common-core
completion~\eqref{eq:completion}, with idempotent $e$.
Put $f=I-e$, assume $\rk f=3$, and suppose
$T_{\rm core}=e(a_0-a_0b_0a_0)e$ is invertible in $eM_n(\C)e$.
Identify $f\C^n$ with $\C^3$ and suppose
$N=fa_0f$, $C=fc_0f$ are nilpotent, $N+C$ is invertible,
$\beta=fb_0f=(N+C)^{-1}$, and $\nu=N\beta C\ne0$ is nilpotent.
For analytic raw paths through $a_0,b_0$, the three eigenvalues of
$T_\varepsilon=a_\varepsilon-a_\varepsilon b_\varepsilon a_\varepsilon$
tending to zero are $O(|\varepsilon|^{2/3})$.
Off-corner perturbations are allowed, and no continuous
primed-inverse path is assumed.
\end{corollary}
\begin{proof}
A small circle separating zero from the limiting core spectrum
defines an analytic Riesz projection $E(\varepsilon)$ with $E(0)=f$.
Set
\[
 W(\varepsilon)=E(\varepsilon)f+(I-E(\varepsilon))(I-f).
\]
Then $W(0)=I$ and $E(\varepsilon)W(\varepsilon)=W(\varepsilon)f$,
so $W$ is locally invertible. On $f\C^n$, the matrix
$\widehat T_\varepsilon=fW(\varepsilon)^{-1}T_\varepsilon
W(\varepsilon)f$ has exactly the selected eigenvalues and
$\widehat T_0=\nu$. With $K=W'(0)$,
$X=\left.\frac{d}{d\varepsilon}a_\varepsilon\right|_0$ and
$Y=\left.\frac{d}{d\varepsilon}b_\varepsilon\right|_0$,
\[
 \widehat T'_0=fT'_0f+[\nu,fKf],\qquad
 fT'_0f=X_f-X_f\beta N-N\beta X_f-NY_fN,
\]
where $X_f=fXf$, $Y_f=fYf$. The second formula uses commutation
of the limiting input data with $f$, and does not assume that the
perturbed data commute with $f$.
Theorem~\ref{thm:mechanism}(4) and trace cyclicity give
$\tr(\nu\widehat T'_0)=0$. Hence the linear characteristic
coefficient of the rank-one square-zero compressed block is
$O(\varepsilon^2)$. Its determinant is also
$O(\varepsilon^2)$, since every $2\times2$ minor of $\nu$
vanishes. Its trace is $O(\varepsilon)$. The cubic rescaling
in Theorem~\ref{thm:bound} now applies.
\end{proof}

\section{Rates, gaps, and the additive regime}\label{sec:scope}

\subsection{Eigenvalue scales and spectral gaps}
The largest small-eigenvalue modulus and the punctured spectral gap
can have different orders. For nonzero leading constants and exponents
$r_1,\ldots,r_k$, their orders are respectively
$|\varepsilon|^{\min r_j}$ and $|\varepsilon|^{\max r_j}$.
Here ``order'' means two-sided positive bounds for small nonzero
parameter on a local branch. The following table records proved cases;
its last two columns give the exponents of these two quantities.
\begin{center}
\small
\begin{tabular}{@{}p{0.47\linewidth}cc@{}}
\hline
Base and path condition & Maximum & Gap\\ \hline
$m=3$, type $[2,1]$, $\kappa_{\rm raw}\ne0$ & $2/3$ & $2/3$\\
Canonical $m=3$, $\alpha=0$ & $2$ & $2$\\
Canonical $m=3$, $\alpha\ne0$ & $1$ & $5/2$\\
$m\ge5$, type $[m-1,1]$, $c_1\kappa_{\rm raw}\ne0$
  & $1/(m-2)$ & $1/2$\\
$m=4$, type $[3,1]$, $c_1\kappa_{\rm raw}\ne0$ & $1/2$ & $1/2$\\
$m=4$, type $[2,2]$, fixed coupling, $\kappa_4\ne0$ & $1/2$ & $1/2$\\
Zero-primary multiplicity two, $\kappa_{\rm raw}\ne0$ & $1$ & $1$\\
The sharp exceptional path~\eqref{eq:two-zero-sharp-cubic}
  & $3/2$ & $3/2$\\
The exceptional affine path of Example~\ref{ex:affine-cubic-gap}
  & $1$ & $3$\\ \hline
\end{tabular}
\end{center}
Every row assumes nilpotent base inputs. The canonical rows also
assume a nilpotent base residue, possibly zero.
Theorems~\ref{thm:generic} and~\ref{thm:canonical} give the
three-dimensional generic and canonical exponent lists
$[2/3,2/3,2/3]$, $[2,2,2]$ and $[1,5/2,5/2]$.
The extremal-type rows use Theorem~\ref{thm:extremal-splitting},
and the type-$[2,2]$ row uses Proposition~\ref{prop:full-cancellation}.
The zero-primary multiplicity-two and sharp $3/2$ rows concern
only the cluster tending to zero in Corollary~\ref{cor:other-rates}.
Its exceptional $O(|\varepsilon|^{3/2})$ bound need not be attained
on every path; equation~\eqref{eq:two-zero-sharp-cubic} is the
specified sharp example. The exceptional affine path in
Example~\ref{ex:affine-cubic-gap} has one permanent zero branch,
excluded from the gap, and two nonzero branches of orders $1$ and $3$.
Corollary~\ref{cor:affine-sharp-gap} shows that this cubic gap is
extremal over the exceptional affine class. In the subclass with no
permanent zero, the maximal gap exponent is $5/2$, attained by the
canonical example in Remark~\ref{rem:census}.
In the nilpotent-preserving family~\eqref{eq:nilpotent-preserving-sharpness},
two roots remain zero and both the maximum and the punctured gap
of the other roots have order $|\varepsilon|^{1/(m-2)}$.
When the small cluster contains nonzero eigenvalues tending to zero,
adjoining a fixed invertible core leaves the full gap determined by that
cluster for sufficiently small nonzero parameter. For an invertible
matrix $T$, spectral mapping gives
\[
 \rho(T^{-1})=\gap(T)^{-1}\le\|T^{-1}\|.
\]
This is an equality for the inverse spectral radius and a lower bound for the inverse norm. A general upper norm estimate requires information beyond the spectral gap.

\begin{example}[A flat smooth gap]\label{ex:flat}
The quadratic complementary law of Proposition~\ref{prop:quadratic}
also gives a smooth path with an independently moving scalar core.
Set $e=E_{11}$ in $M_3(\C)$ and hold $b^{-1}=2e+E_{23}+E_{32}$ fixed. For real $\delta\ne0$, put $\eta(\delta)=\exp(-1/\delta^2)$ and set $\eta(0)=0$. The smooth inputs
\[
 a_\delta=(1+\delta)e+E_{23}+\eta(\delta)(E_{22}+E_{33}),\qquad
 c_\delta=b^{-1}-a_\delta
\]
satisfy, for small $\delta$,
\[
 \|a_\delta-a_0\|=\|c_\delta-c_0\|
 =\max\{|\delta|,\eta(\delta)\}=|\delta|.
\]
Their output is
\[
 t_\delta=\frac{1-\delta^2}{2}e-\eta(\delta)^2(E_{23}+E_{32}).
\]
For sufficiently small nonzero $\delta$ it is invertible and
$\gap(t_\delta)=\exp(-2/\delta^2)$. Thus the gap closes faster than every positive power of the raw input displacement. The path is smooth but not analytic at zero. Remaining finite-order rate questions require a specified dimension and analytic, algebraic, or otherwise restricted path or direction class.
\end{example}

In the three-dimensional exceptional affine class,
Proposition~\ref{prop:affine-coefficient-criterion} computes the exponents;
Problem~\ref{prob:rate} concerns the remaining geometric classification.

\subsection{Inverse selection and continuity}
The following is the finite-dimensional form of
Badea--Mbekhta~\cite[Theorem~3.2]{BadeaMbekhta97}. Their operator
theorem states a supremum; the matrix argument below also proves
attainment by a reflexive inverse.

\begin{proposition}[A constrained spectral-gap formula]\label{prop:output-primary-gap}
Let $T\in M_m(\C)$, $m\ge1$, and let $X_*(T)$ be the sum of its
primary spaces for nonzero eigenvalues. Then
\[
 \gap(T)=
 \max_{\substack{S\in M_m(\C),\ TST=T\\SX_*(T)\subseteq X_*(T)}}
       \frac1{\rho(S)},\qquad 1/0=+\infty.
\]
The maximum is attained by an $S$ also satisfying $STS=S$.
\end{proposition}
\begin{proof}
If $X_*(T)\ne\{0\}$, write $T=M\oplus J$ in primary
coordinates, where $M$ is invertible and $J$ is nilpotent;
the nilpotent summand may be absent. Every allowed $S$ has the form
\[
 S=\begin{pmatrix}M^{-1}&B\\0&D\end{pmatrix},
 \qquad BJ=0,\qquad JDJ=J.
\]
Consequently $\rho(S)\ge\rho(M^{-1})=\gap(T)^{-1}$.
In Jordan coordinates for $J$, let $L$ be the direct sum of the
reverse shifts on its Jordan blocks. Then $L$ is nilpotent,
$JLJ=J$ and $LJL=L$, so $M^{-1}\oplus L$ is a reflexive
inverse attaining the bound. If $T$ is nilpotent, the same reverse
shifts give a nilpotent reflexive inverse and attain $+\infty$.
\end{proof}

\emph{Why output-primary invariance is needed.}
The invariance condition is essential even for a residue of two
nilpotent inputs. In the scalar specialization of~\eqref{eq:Q-residue},
fix $T=R_2=\tfrac12E_{12}+E_{22}$ and put
\[
 S_\delta=(2-2\delta)E_{21}+\delta E_{22}\qquad(\delta\in\C).
\]
Direct multiplication gives $TS_\delta T=T$,
$S_\delta TS_\delta=S_\delta$ and $S_\delta^2=\delta S_\delta$;
also $\det(zI-S_\delta)=z^2(z-\delta)$. Thus
$\gap(T)=1$, whereas $1/\rho(S_\delta)=1/|\delta|$ is
unbounded for nonzero $\delta\to0$. Here
$X_*(T)=\operatorname{span}(e_1+2e_2)$ and
$S_\delta(e_1+2e_2)=2e_2\notin X_*(T)$.
The realizing nilpotent inputs have common core zero, whose
preservation imposes no condition. It is the output's nonzero
primary space that must be preserved in the gap formula.

The raw-input splitting paths considered above need not lift continuously
to the prescribed canonical inverse coordinates within the finite-Drazin
Hua class, as Example~\ref{ex:twothirds} shows. In complex Banach algebras,
\cite[Theorem~8.3]{HuaI} characterizes continuity of the core and
candidate along convergent raw pairs whose terms and limit are
full configurations: intertwining holds and the input remainders
are quasinilpotent. Its criterion uses uniform input spectral gaps.
Theorem~8.4 of that paper adds an equivalent output-gap criterion
on the restricted class, where the residues at every term and at
the limit are also quasinilpotent. Proposition~9.6 and Remark~9.7
concern local variation of augmented outer-inverse data.
The splitting results here concern small output eigenvalues along
specified raw paths and assert neither continuity of the inverse
coordinates nor a uniform inverse-norm bound.

\emph{Fixed reflexive inversion and discontinuous Drazin inversion.}
For the scalar family $R_x=\tfrac12E_{12}+\tfrac x2E_{22}$,
the identity $R_x^2=(x/2)R_x$ gives
\[
 R_x^{\Dz}=2x^{-2}E_{12}+2x^{-1}E_{22}\quad(x\ne0),
 \qquad R_0^{\Dz}=0.
\]
For $x\ne0$,
\[
 \gap(R_x)=|x|/2,\qquad \rho(R_x^{\Dz})=2/|x|,\qquad
 \|R_x^{\Dz}\|=\frac{2\sqrt{1+|x|^2}}{|x|^2}.
\]
The fixed reflexive inverse $G=2E_{21}$ remains bounded while
the Drazin inverse diverges, and its norm grows faster than its
spectral radius. Both realizing inputs remain nilpotent with
Drazin inverses zero, so their common core and the prescribed
core candidate are constantly zero. The candidate succeeds at
$x=0$ and fails for $x\ne0$, although the actual output Drazin
inverse exists for every $x$. For $x\ne0$ the residue is not
quasinilpotent, so this path leaves the restricted class required
for the output-gap continuity criterion just cited.

\subsection{The additive comparison}
The requirement that the two remainders have invertible sum also identifies the relevant additive regime. For Banach generalized Drazin invertible elements $x,y$, the compatibility conditions in Castro-Gonz\'alez and Koliha~\cite[equation (3.1)]{CGK04} are
\[
 x^\pi y=y,\qquad xy^\pi=x,\qquad y^\pi xyx^\pi=0,
 \qquad x^\pi=1-xx^{\gd},\quad y^\pi=1-yy^{\gd}.
\]
For quasinilpotent inputs, $x^\pi=y^\pi=1$, and these conditions reduce to $xy=0$. For every nonzero complex $\lambda$,
\[
 (\lambda-x)(\lambda-y)=\lambda(\lambda-x-y).
\]
The factors on the left are invertible, so $x+y$ is quasinilpotent and cannot be invertible in a nonzero Banach algebra. This is the quasinilpotent specialization used in the first case of their Theorem~3.3, through Lemma~2.4. The comparison concerns those actual conditions, not all additive generalized-Drazin results.

In our complementary corner, an arbitrary pair of quasinilpotent inputs with invertible sum can be completed by taking $a=e+N$, $c=e+C$, $a'=c'=e$ and $b=\tfrac12e+(N+C)^{-1}_{f\A f}$ in a complex unital Banach algebra. The equations of Section~\ref{sec:corners} verify the completion. Conversely, common-core coupling makes that sum invertible in its corner. Thus input-pair completion is settled. Quasinilpotence of $\nu$
remains the separate condition for success of the prescribed core inverse.
Theorem~\ref{thm:complete-realization} settles the matrix residue set over every infinite field
in every dimension $m\ge2$. The coefficient construction
and Banach consequences in Theorem~\ref{thm:realization} and
Example~\ref{ex:Banach} retain their separate scope.

\section{Remaining questions}\label{sec:problems}
The infinite-field residue set, the generic extremal laws, the
universal full-cancellation criterion and the finite exceptional
affine coefficient test are settled here. Questions remain about
the geometry of exceptional directions and of realizing pairs.
The finite-field obstruction in Remark~\ref{rem:realization-F2}
does not provide a classification over other finite fields.

\begin{problem}[Geometry of exceptional affine directions]\label{prob:rate}
Fix nilpotent $N,C\in M_3(\C)$ with invertible sum and nonzero
nilpotent residue. For a specified complex linear space
$L\subseteq M_3(\C)$ with $\kappa|_L=0$, consider
\[
 A_\varepsilon=N+\varepsilon X,\qquad
 C_\varepsilon=C-\varepsilon X,\qquad X\in L.
\]
Proposition~\ref{prop:affine-coefficient-criterion} computes the
exponents for each $X$ and the generic exponent list on $L$.
Determine which valuation and exponent loci occur, their dimensions
and closure relations, and how these depend on the base and on $L$.
The normal form in Lemma~\ref{lem:normal-form} places $L$ in one
of the two exceptional hyperplanes; the geometry of its further
coefficient loci is not classified by the finite test.
For general analytic paths, a corresponding question must also
specify the allowed higher input coefficients.
\end{problem}

\begin{problem}[Compatible singular lifting]\label{prob:monodromy}
The unconstrained roots of nonsingular $K=I-4bt$ are classified in
\cite[Section~1.6, Theorems~1.24 and~1.26]{HighamBook}.
Beyond the simple nonzero-spectrum base of
Proposition~\ref{prop:monodromy}, determine which roots and analytic
selections are compatible with a specified spectral-configuration
stratum (Remark~\ref{rem:raw}) when eigenvalues of $K$ collide or
reach zero. The generic monodromy on that base is settled.
Singularity of the target $t$ means $1\in\sigma(K)$; singularity
of the square-root target $K$ means $1/4\in\sigma(bt)$.
These are different conditions.
\end{problem}

\begin{problem}[Geometry of nonnilpotent residues]\label{prob:regime}
For a prescribed nonnilpotent residue, determine which further
kernel/image incidences and ordered annihilation patterns are
possible or forced among its nilpotent realizing pairs, beyond
Lemma~\ref{lem:rank} and the rank-one zero-primary conclusion of
Section~\ref{sec:mechanism}. Theorem~\ref{thm:complete-realization}
settles which matrices occur. It does not classify all pairs
realizing a fixed target.
\end{problem}

\appendix
\section{Square roots of the raw Hua equation}\label{sec:cover}
This appendix concerns $t=a-aba$ with $b,t$ prescribed,
before imposing compatibility conditions on generalized-inverse data.
The root counts, critical-point criterion, and covering statement
use their stated finite-dimensional hypotheses; the additional
spectral-intertwining condition is discussed separately.

\begin{proposition}[The algebraic square-root correspondence]\label{prop:cover}
Let $\R$ be a unital ring in which $2$ is invertible, let $b\in\R^\times$, and prescribe $t\in\R$. The raw equation $t=a-aba$ is equivalent, through
\[
y=ba,\qquad X=1-2y,
\]
to $X^2=1-4bt$. Its inverse parametrization is $a=\tfrac12b^{-1}(1-X)$. The exchange $a\leftrightarrow c=b^{-1}-a$ is $X\leftrightarrow-X$. Every solution $y$ commutes with $bt$.

For $M_n(\C)$, if $bt$ has $n$ distinct eigenvalues, none equal to $1/4$, there are exactly $2^n$ raw solutions. At $t=0$ the raw solutions are $b^{-1}$ times the idempotents, so the fiber is infinite for $n\ge2$. These counts do not impose intertwining or spectral-configuration hypotheses.
\end{proposition}
\begin{proof}
Multiplication on the left by $b$ gives $bt=y-y^2$. Completing the square proves the equivalence and its inverse. The exchange replaces $y$ by $1-y$. Commutation follows because $bt$ is a polynomial in $y$. In the stated matrix case, a square root commutes with $K=I-4bt$. This matrix has distinct nonzero eigenvalues, so its commutant is diagonal in an eigenbasis and each eigenline has two square-root choices. At $t=0$ the equation is $y^2=y$. For $n\ge2$, the idempotents $E_{11}+sE_{12}$, $s\in\C$, already give infinitely many distinct solutions after left multiplication by $b^{-1}$.
\end{proof}

\begin{remark}[Raw solutions and spectral configurations]\label{rem:raw}
The finite count and the possible infinite matrix-root fibers are standard phenomena; see Higham~\cite[Sections~1.4 and~1.6, Theorems~1.24 and~1.26]{HighamBook} for roots of nonsingular matrices. For the present matrix comparison, the spectral-intertwining condition is $pb=br$ with $p=aa^{\Dz}$ and $r=cc^{\Dz}$. If $t$ is invertible in the matrix case, $t=abc$ forces $a,c$ invertible by determinants, so every raw solution is classical. Without this extra condition, even the finite raw count need not count spectral configurations. For $b=I_2$ and $t=\operatorname{diag}(0,-2)$ there are four raw solutions. Indeed $K=I-4bt=\operatorname{diag}(1,9)$, so all roots and raw solutions are diagonal in these coordinates. The first coordinate of $a$ is $0$ or $1$, and that of $c$ is its complement, so their Drazin core projections differ in that coordinate. None satisfies the spectral intertwining condition for $b=I_2$.

On the simple nonzero-spectrum base of Proposition~\ref{prop:monodromy},
this gives a general dichotomy. If $t$ is invertible, every raw solution
is classical, as above. If $t$ is singular, $K=I-4bt$ has eigenvalue $1$
on exactly one eigenline. Since $X^2=K$ implies $XK=KX$, each root is
diagonal in an eigenbasis of $K$: it chooses $+1$ or $-1$ on that line
and has no other eigenvalue $\pm1$. Hence precisely one of
$a=\tfrac12b^{-1}(I-X)$ and $c=\tfrac12b^{-1}(I+X)$ is singular, and the
other is invertible. The invertible input has canonical Drazin core $I$.
If $p=I$, the relation $pb=br$ and invertibility of $b$ force
$r=I$; if $r=I$, they force $p=I$. Either conclusion contradicts
singularity of the other input. Thus no raw solution here satisfies
spectral intertwining.
No commutation of $b$ with $X$ is used. This dichotomy is confined to the
specified simple-spectrum base.
\end{remark}

\begin{proposition}[Critical points and analytic selection]\label{prop:critical}
For fixed $b\in GL_n(\C)$, the derivative of $a\mapsto a-aba$ is invertible at $a$ if and only if
\[
\sigma(X)\cap(-\sigma(X))=\varnothing,\qquad X=I-2ba.
\]
Thus singularity of $a-c$ is not the entire critical locus. Pointwise solvability of the raw equation along a path does not imply a holomorphic choice of inputs.
\end{proposition}
\begin{proof}
For a direction $H$ put $Z=bH$. Left multiplication by $b$ transforms the derivative into
\[
Z-Zy-yZ=\tfrac12(XZ+ZX),\qquad y=ba.
\]
The operator $L_X:Z\mapsto XZ+ZX$ has matrix $I\otimes X+X^T\otimes I$ under column vectorization. Triangularizing the two tensor factors shows that its eigenvalues are the sums $\lambda_i(X)+\lambda_j(X)$, counted with multiplicity. This proves the criterion, the finite-dimensional Sylvester separation condition; see~\cite[Section~2]{BhatiaRosenthal}. When the derivative is invertible, the complex analytic inverse-function theorem applied to the polynomial map $a\mapsto a-aba$ supplies a local analytic inverse near that point. This local existence argument is separate from the failure of analytic selection at the critical point below.

For a concrete critical point, take $b=I_2$ and
$a_s=\bigl(\begin{smallmatrix}2&s\\0&-1\end{smallmatrix}\bigr)$. Then $a_s-a_s^2=-2I_2$ for every $s$, while $a_s-c_s=\bigl(\begin{smallmatrix}3&2s\\0&-3\end{smallmatrix}\bigr)$ is invertible. The corresponding $X$ has eigenvalues $-3,3$, so $L_X$ is singular. This is an infinite classical fiber.
For analytic selection take scalars $b=1$ and $t(z)=(1-z)/4$. The square-root equation is $X(z)^2=z$. Pointwise solutions exist, but a holomorphic solution at zero would have a square with even vanishing order, a contradiction. Here $t(z)\ne0$ locally, so even pointwise classical solvability does not repair the lifting problem.
\end{proof}

\begin{proposition}[The generic square-root cover]\label{prop:monodromy}
On the open set of raw base pairs $(b,t)$ for which $b$ is invertible and $K=I-4bt$ has distinct nonzero eigenvalues, the raw solution map is a $2^n$-sheeted unbranched cover. Its monodromy on the set of sign choices is the signed-permutation action of $(\mathbb Z/2\mathbb Z)^n\rtimes S_n$. This is a statement on that specified open set, not a global fibration or a critical-locus description.
\end{proposition}
\begin{proof}
The change of variables $(b,t)\leftrightarrow(b,K)$ identifies the base with $GL_n(\C)$ times the simple nonzero-spectrum locus. Locally the simple eigenvalues $\mu_i$ are analytic by the implicit-function theorem for the characteristic polynomial. Their spectral projections are the analytic matrices
\[
 E_i=\prod_{j\ne i}\frac{K-\mu_j I}{\mu_i-\mu_j}.
\]
Choosing local square-root branches for the nonzero $\mu_i$, every square root has the form
\[
X=\sum_{i=1}^n\epsilon_i\sqrt{\mu_i}\,E_i,\qquad \epsilon_i\in\{1,-1\}.
\]
The $2^n$ expressions are distinct local analytic branches. Continuation around any loop permutes the eigenlines and may change their square-root signs, so the monodromy is contained in the signed-permutation group. A diagonal loop moving one eigenvalue once around zero while avoiding the other fixed eigenvalues changes that sign alone; take $b$ fixed and recover $t$ from $K$. Paths interchanging two eigenvalues while avoiding zero, combined with a path of basis changes ending at the corresponding permutation matrix, give loops whose eigenline permutations are transpositions. Their possible additional sign changes can be canceled by the first loops. These operations generate the full signed-permutation action.
\end{proof}

\section*{Acknowledgments}
The author would like to acknowledge the sincere support provided by the library facilities of the Harish-Chandra Research Institute (HRI), Prayagraj (Allahabad), India, during the author's five-year visiting period at the institute. This was the time when principal ideas pertaining to this program were conceived. This support was greatly facilitated by the then Director, Prof. Pinaki Majumdar, to whom the author is infinitely grateful. The author is also thankful to Prof. Ilya Spitkovsky for his academic support.

\end{document}